\documentclass{article}

\usepackage{comment}
\usepackage{amsfonts}
\usepackage{color}
\usepackage{amssymb}
\usepackage{verbatim}
\usepackage{amsmath}
\usepackage{esint}
\usepackage{amscd}
\usepackage{euscript}
\usepackage{mathtools}
\usepackage[all]{xy} 
\usepackage{graphicx}
\usepackage{hyperref}
\usepackage{graphicx}
\usepackage{rotating}
\usepackage[usenames,dvipsnames]{xcolor}
\usepackage{datetime}
\usepackage{mathtools}
\usepackage{tcolorbox}
\usepackage{cancel}

\newtheorem{Theorem}{Theorem}
\newtheorem{Corollary}[Theorem]{Corollary}
\newtheorem{Example}[Theorem]{Example}
\newtheorem{Definition}[Theorem]{Definition}
\newtheorem{remark}[Theorem]{Remark}

\newtheorem{Lemma}[Theorem]{Lemma}
\newtheorem{Proposition}[Theorem]{Proposition}

\newtheorem{Fundamental Theorem}{Fundamental Theorem}

\newenvironment{Proof}[1][Proof]{\textbf{#1.} }{\ \rule{0.5em}{0.5em}}

\newenvironment{Notation}[1][Notation]{\textbf{#1.} }

\newenvironment{Idea}[1][Idea]{\textbf{#1.} }

\newenvironment{question}[1][Question]{\textbf{#1.} }

\usepackage[all]{xy}

\def \C {\mathfrak{C}}
\def \B {\mathcal{B}}

\def \bD {\underline{\Delta}}
\def \bS {\underline{S}}
\def \bba {{\underline{\underline{a}}}}
\def \ba {{\underline{a}}}

\def \btn {\hat{\tn}}

\def \bad {\underline{\tad}}
\def \u {\underline}
\def \ad {\triangleright_{ad}}
\def \bad {\blacktriangleright_{bad}}

\def \q {\quad}
\def \cb {\textcolor{black}}

\def \we {{\textrm{ where }}}

\def \an {\textrm{ and }}
\def \dis {\displaystyle}
\def \t {\triangleright}
\def \tr {\triangleright_\rho}
\def \tn {\otimes}
\def \id {{\rm id}}

\def \r {{{\rm REqual}}}

\def \lk {{{\rm LKer}}}
\def \rk {{{\rm RKer}}}
\def \ck {{{\rm CKer}}}

\def \D {\Delta}
\def \e {\epsilon}
\def \z  {\zeta}
\def \k {\kappa}
\def \tr {\triangleright}
\def \tl {\triangleleft}
\def \d {\partial}
\def \F {\mathcal{F}}

\def \n {\noindent}

\def \yd {\mathcal{YD}_{\C}(H)}

\def \ra {\xrightarrow}

\allowdisplaybreaks

\begin{document}
	
\title{Braided Hopf Crossed Modules Through Simplicial Structures}

\author{Kadir Emir \\ emir@math.muni.cz \and Jan Paseka \\ paseka@math.muni.cz \and  {\small Department of Mathematics and Statistics, Masaryk University, Brno, Czech Republic.}}

%\date{\today \quad \currenttime}

\maketitle

\begin{abstract}
	Any simplicial Hopf algebra involves $2n$ different projections between the Hopf algebras $H_n,H_{n-1}$ for each $n \geq 1$. The word projection, here meaning a tuple $\d \colon H_{n} \to H_{n-1}$ and $i \colon H_{n-1} \to H_{n}$ of Hopf algebra morphisms, such that $\d \, i = \id$. Given a Hopf algebra projection $(\d \colon I \to H,i)$ in a braided monoidal category $\C$, one can obtain a new Hopf algebra structure living in the category of Yetter-Drinfeld modules over $H$, due to Radford's theorem. The underlying set of this Hopf algebra is obtained by an equalizer which only defines a sub-algebra (not a sub-coalgebra) of $I$ in $\C$. In fact, this is a braided Hopf algebra since the category of Yetter-Drinfeld modules over a Hopf algebra with an invertible antipode is braided monoidal. To apply Radford's theorem in a simplicial Hopf algebra successively, we require some extra functorial properties of Yetter-Drinfeld modules. Furthermore, this allows us to model Majid's braided Hopf crossed module notion from the perspective of a simplicial structure.
\end{abstract}

\medskip


\noindent{\bf 2010 AMS Classification}
{16T05, 18D05, 18D10, 18G30.}

{\small \tableofcontents}

\section{Introduction}

A Hopf algebra \cite{SW1} can be considered as an abstraction of the group algebra (of a group) and the universal enveloping algebra (of a Lie algebra). In other words, considering the group-like and primitive elements in a Hopf algebra, we have the functors: 
\begin{align}\label{functors1}
\begin{split}
\mathrm{( \ )_{gl}^{*} \colon \{Hopf \ Algebras\}} & \to \mathrm{ \{Groups\} } \, ,\\
\mathrm{Prim \colon \{Hopf \ Algebras\}} & \to \mathrm{\{Lie \ Algebras\} } \, .
\end{split}
\end{align}
 There exists a wide variety of variations of Hopf algebras by  relaxing its properties or adding some extra structures; for instance, quasi-Hopf algebras \cite{zbMATH07003760}, quasi-triangular Hopf algebras \cite{SM1} and quantum groups \cite{SM1, zbMATH05015420}, etc. The original definition of Hopf algebras is given in the category of vector spaces. However, we know that the category of vector spaces is a trivial example of (symmetric) braided monoidal categories. As a consequence of this categorical relationship, a more general definition of Hopf algebras in a braided monoidal category is given in \cite{SM3}.

\medskip

Suppose that, we have a Hopf algebra projection $(\d \colon I \to H,i)$ in the category of vector spaces (i.e. there exists a Hopf algebra map $i \colon I \to H$ such that $d \, i = \id_{H}$). One can obtain a Yetter-Drinfeld module \cite{zbMATH04171177} over $H$ where the action and coaction is defined by taking the Hopf algebra morphisms $\d$ and $i$ into account. The category of such Yetter-Drinfeld modules will be denoted by $\mathcal{YD}(H)$. Moreover, there exists a Hopf algebra structure $\B$ (subset of $I$) living in $\mathcal{YD}(H)$. The underlying set of the Hopf algebra $\B$ is categorically given by an equalizer which is only a sub-algebra of $I$. In fact, this is a braided Hopf algebra since the category of Yetter-Drinfeld modules over a Hopf algebra with an invertible antipode is braided monoidal. On the other hand, considering the underlying set of $\B$, we have a Hopf algebra structure $\B \btn H$ in the category of vector spaces (not in $\mathcal{YD}(H)$) comes equipped with the smash product and smash coproduct. Furthermore, we have the isomorphism of Hopf algebras $I \cong \B \btn I$. This is called ``Radford's theorem" which is introduced in \cite{Rad1} and extended in \cite{zbMATH00487164}. Inspired by the categorical relationship between vector spaces and braided monoidal categories, it is proven in \cite{zbMATH01036413} that, Radford's Theorem is also valid in any braided monoidal category $\C$ with equalizers. 

\medskip

The notion of strict 2-group is a strict vertical categorification of that of group. Equivalently, a strict 2-group is the internal category in the category of groups. The essential example of a strict 2-group is called crossed module  \cite{WCH2}. In other words, they can be seen as a way to encode a strict group \cite{zbMATH03521203}. For more details on crossed modules especially from the topological, algebraic and geometric point of view, we refer \cite{zbMATH06825034,zbMATH05815821,zbMATH06520714}. Moreover, the notion and some applications of crossed modules are extended to various algebraic structures such as monoids \cite{zbMATH07140994,Bohm2}, Lie algebras \cite{E1}, groupoids \cite{zbMATH01944940}, groups with operations \cite{zbMATH03957392}, modified categories of interest \cite{zbMATH06456664}, shelves (via generalized Yetter-Drinfeld modules) \cite{zbMATH06812675}, etc. More recently, motivated from the group theoretical analogy, Majid defines strict quantum 2-groups in \cite{SM2}. Through encoding that of strict quantum 2-groups, his paper also includes two different approaches to crossed module notion for Hopf algebras, as follows:

\newpage

\begin{itemize}
	\item Hopf crossed modules:
	\begin{itemize}
		\item They are given by a Hopf algebra morphism in the category of vector spaces.
		\item They unify the Lie algebra crossed modules and group crossed modules with the following functors (between the categories of crossed modules of Lie algebras, Hopf algebras and groups):
			\begin{align*}
		\xymatrix@R=40pt@C=30pt{ \mathrm{\{XLie\}}  &
			\mathrm{\{XHopf\}} \ar[r]^{( \ )_{gl}^{*}} \ar[l]_{\mathrm{Prim}} &
			\mathrm{\{XGrp\}} \, .
		}
		\end{align*}
		which are induced from \eqref{functors1}. For more details about the functorial relationships between various algebraic categories of crossed modules, see \cite{JFM1}.
		\item Categorically, they are equivalent to the certain case of simplicial Hopf algebras given in \cite{1905.09620}, if we impose the cocommutativity condition.
		\item We also know that the category of cocommutative Hopf algebras over a field of characteristic zero is semi-abelian. From this point of view, it is proven in \cite{zbMATH07075977} that Majid's crossed module definition is coherent in the sense of internal crossed modules in semi-abelian categories \cite{zbMATH02096926}.
	\end{itemize}
	\item{{Braided Hopf crossed modules:}}
	\begin{itemize}
		\item They are given by a twisted Hopf algebra map: briefly, it is not a Hopf algebra morphism in a certain category.
		\item The codomain of braided Hopf crossed modules are still vector spaces, while the domains are now Yetter-Drinfeld modules.
		\item The main motivation to introduce this notion is to generalize the characterization of Hopf crossed modules.
	\end{itemize}
\end{itemize}

\medskip

When we forget the cocommutativity case, we notice that Majid's first crossed module notation does not make a contact with a simplicial structure in the sense of \cite{1905.09620}. Because, the structure $\mathcal{B}$ we mentioned in the Radford's theorem is: 
\begin{itemize}
	\item a Hopf algebra living in the base category (vector spaces), if the Hopf algebras are cocommutative,
	\item a Hopf algebra living in the category of Yetter-Drinfeld modules, in general case (without the cocommutativity condition).
\end{itemize}

\medskip

So, it gives us an idea that, Majid's braided Hopf crossed module notion would be the correct one to model crossed modules in the sense of simplicial structures and also to discover its other well-known group theoretical properties in the category of Hopf algebras. For this aim, in this paper, we first recall Radford's theorem in a braided monoidal category $\C$. Afterwards, we see that a Hopf algebra living in the braided monoidal category $\mathcal{YD}_{\C}(H)$ can also be converted into a Hopf algebra living in the braided monoidal category $\mathcal{YD}_{\C}(I)$, for a given Hopf algebra projection $(\d \colon I \to H,i)$. This property will be used to obtain new projections in induced categories in a simplicial Hopf algebra. Consequently, it gives rise to a connection between braided Hopf crossed modules and simplicial Hopf algebras. 

\medskip

\section*{Acknowledgment}

The first author acknowledges the support by the Austrian Science Fund (FWF): project I 4579-N and the Czech Science Foundation (GA\v CR): project 20-09869L. He also acknowledges the support
by the project Mathematical structures 9 (MUNI/A/0885/2019). The second author was supported by the project Group Techniques and Quantum Information (MUNI/G/1211/2017) by Masaryk University Grant Agency (GAMU).

\medskip

\section{Preliminaries}

\subsection{Braided monoidal categories}

A braided monoidal category \cite{JS2} (or braided tensor category) is a monoidal category $\C$ equipped with a braiding map (natural isomorphism) $ \mathcal{R}_{A,B} \colon A \tn B \to B \tn A$ for each objects of $\C$, satisfying the following hexagon diagrams:
\begin{align*}
\xymatrix@R=10pt@C=15pt{
	& A \tn (B \tn C) \ar[rr]^-{\mathcal{R}_{A,B \tn C}} & & (B \tn C) \tn A  \ar[dr]^-{\alpha_{B,C,A}}& 
	\\ (A \tn B) \tn C \ar[ur]^-{\alpha_{A,B,C}}  \ar[dr]_-{ \mathcal{R}_{A,B} \tn \id} & & & & B \tn (C \tn A) 
	\\ & (B \tn A) \tn C \ar[rr]^-{\alpha_{B,A,C}} & & B \tn (A \tn C)  \ar[ur]_-{\id \tn \mathcal{R}_{C,A}} &  } \\ \\
\xymatrix@R=10pt@C=15pt{ 
	& (A \tn B) \tn C \ar[rr]^-{\mathcal{R}_{A\tn B, C}} & & C \tn (A \tn B)  \ar[dr]^-{\alpha^{-1}_{C,A,B}}& 
\\ A \tn (B \tn C) \ar[ur]^-{\alpha^{-1}_{A,B,C}}  \ar[dr]_-{\id \tn {\mathcal{R}_{B , C}}} & & & & (C \tn A) \tn B 
\\ & A \tn (C \tn B) \ar[rr]^-{\alpha^{-1}_{A,C,B}} & & (A \tn C) \tn B  \ar[ur]_-{\mathcal{R}_{A,C} \tn \id} &  } 
\end{align*}
where $\alpha_{A,B,C} \colon (A \tn B) \tn C \to A \tn (B \tn C)$ denotes the associator.

\medskip

	Any symmetric monoidal category is a braided monoidal category. Therefore, the category of vector spaces is a braided monoidal category with a trivial braiding $ \mathcal{R}_{V,W} \colon v \tn w \in V \tn W \mapsto w \tn v  \in W \tn V$. As another example with a non-trivial braiding, suppose that $G$ is a group. Then, the category of crossed $G$-sets forms a braided monoidal category where the braiding is given by using the group action.

\subsection{Hopf algebraic conventions}

In this subsection, we recall some Hopf algebraic notions in a braided monoidal category $\C$ from \cite{SM3, zbMATH06738086} which generalizes the ordinary Hopf algebras \cite{SM1} (i.e. Hopf algebras in the category of vector spaces). The most important examples are the braided Hopf algebras living in the category of Yetter-Drinfeld modules (\S \ref{braided-in-yd}), and the Nichols algebras \cite{zbMATH01824088}. 

%\subsubsection{Hopf algebras}

\bigskip

Let $H=(H,\nabla,\eta,\Delta,\epsilon,S)$ be a Hopf algebra in $\C$. That means:
\begin{itemize}
	\item $(H,\nabla,\eta)$ is a unital associative $\k$-algebra. Thus
	\begin{align*}
	\xymatrix@R=20pt@C=30pt{ H \tn H \tn H \ar[r]^-{\mathrm{\nabla \tn \id}} \ar[d]_{\mathrm{\id \tn \nabla}} & H \tn H \ar[d]^{\nabla} \\ H \tn H \ar[r]^-{\nabla} & H }  \quad \quad \quad \quad
	\xymatrix@R=20pt@C=30pt{ \k \tn H \ar[r]^-{\mathrm{\eta \tn \id}} \ar[dr]_{\cong} & H \tn H \ar[d]_-{\nabla} & H \tn \k \ar[dl]^{\cong} \ar[l]_-{\mathrm{\id \tn \eta}}  \\ & H &  } 
	\end{align*}
%	\begin{itemize}
%		\item $\nabla\colon H \tn H \to H$ is an associative product, making $H$ into an associative algebra. In short the product in $H$ induces a map $\nabla\colon H\tn H \to H$, where $x\tn y\mapsto xy$.  
%		\item $\eta\colon \kappa \to H$ is an algebra map endowing $H$ with a unit.  In short $\eta \colon \lambda \in  \C \mapsto \lambda 1_H \in H.$ (Here $1_H$ is the identity element of $H$.)
%	\end{itemize}
	\item $(H,\D,\e)$ is a counital coassociative $\k$-coalgebra. Thus 
	\begin{align*}
	\xymatrix@R=20pt@C=30pt{ H \ar[r]^-{\D} \ar[d]_{\D} & H \tn H \ar[d]^{\D \tn \id} \\ H \tn H \ar[r]^-{\D \tn \id} & H \tn H \tn H }  \quad \quad \quad \quad
	\xymatrix@R=20pt@C=30pt{ \k \tn H \ar@{<-}[r]^-{\mathrm{\e \tn \id}} \ar[dr]_{\cong} & H \tn H  \ar[r]^-{\mathrm{\id \tn \e}} & H \tn \k \ar[dl]^{\cong}   \\ & H \ar[u]^-{\D} &  } 
	\end{align*}
%	\begin{itemize}
%		\item $\D\colon H \to H \tn H$ is a coassociative coproduct. We use Sweedler's notation \cite{SW1} for coproduct:
%		\begin{align*} 
%		\Delta (x) = \underset{(x)}{\sum } \,\, x' \otimes x'', \textrm{ where } x \in H.
%		\end{align*}
%		\item $\epsilon \colon H \to \kappa$ is the counit. So, for all $x \in H$, we have: $$\displaystyle x=\sum_{(x)} \epsilon(x')x''=\sum_{(x)}
%		x'\epsilon(x'').$$
%	\end{itemize}
	\item  $(H,\nabla,\eta,\Delta,\epsilon)$ is a bialgebra; i.e. $\eta$, $\nabla$ are coalgebra morphisms, and $\epsilon$, $\Delta$ are algebra morphisms. Thus
	\begin{align}\label{compa}
	\begin{split}
	& \xymatrix@R=20pt@C=50pt{ H \tn H \ar[r]^-{\D \tn \D} \ar[d]_{\mathrm{\nabla}} & H \tn H \tn H \tn H \ar[r]^-{\id \tn \mathcal{R} \tn \id} & H \tn H \tn H \tn H \ar[d]^-{\nabla \tn \nabla}   
		\\  H \ar[rr]^-{\D }  && H \tn H } 
	\\  & \quad \quad 
	\xymatrix@R=20pt@C=20pt{ \k \ar[r]^-{\mathrm{\eta}} \ar[dr]_{\eta \tn \eta} & H  \ar[d]^{\D} \\ & H \tn H  }  \quad \quad \quad 
	\xymatrix@R=20pt@C=20pt{ \k   & H \ar[l]_-{\mathrm{\e}}  \\ & H \tn H  \ar[u]_{\nabla} \ar[ul]^{\e \tn \e} }  \quad \quad \quad
	\xymatrix@R=20pt@C=20pt{ \k \ar[r]^-{\mathrm{\eta}} \ar@{=}[dr] & H  \ar[d]^{\e} \\ & \k  } 
	\end{split}
	\end{align}
%	where $ \mathcal{R}_{A,B} \colon A \tn B \to B \tn A$ denotes the braiding in $\C$.
	%(In order to prove $a)$ and $b)$ it suffices proving that $a)$ or $b)$ holds.) 
	%\item If $(H,\mu,\eta,\Delta,\epsilon)$ is a bialgebra then the space of linear maps $\hom_\kappa(H,H)$ is itself a unital associative algebra under the convolution product $(f*g)(x)=\sum_{(x)}f(x')g(x'')$. The unit of the convolution algebra is $\eta\circ \epsilon$. And 
	\item There exists an (inverse-like) anti-homomorphism $S\colon H \to H$ at the level of algebra and coalgebra,  satisfying:
	\begin{align}\label{antipode}
	\xymatrix@R=20pt@C=10pt{ 
		& H \tn H \ar[rr]^-{\id \tn S} & & H \tn H  \ar[dr]^-{\nabla}& 
		\\ H \ar[ur]^-{\D}  \ar[dr]_-{\D} \ar[rr]^-{\e} & & \k \ar[rr]^-{\eta} & & H 
		\\ & H \tn H \ar[rr]^-{S \tn \id} & & H \tn H \ar[ur]_-{\nabla} &  } 
	\end{align}
	which is called an antipode.
	%$$\sum_{(x)} S(x')x''=\sum_{(x)} x'S(x'')=\epsilon(x)1_H .$$
\end{itemize}

\medskip

%\begin{Example}
%The most important examples are the braided Hopf algebras in the category of Yetter-Drinfeld modules (\S \ref{braided-in-yd}), and the Nichols algebras \cite{zbMATH01824088}. 
%Here we skipped the well-known definition of a braided monoidal category and refer \cite{JS1} for more details. %Basically, a braided monoidal category $\C$ (or braided tensor category \cite{JS2}) is a monoidal category \cite{selinger} equipped with a braiding map (isomorphism) $ \mathcal{R}_{A,B} \colon A \tn B \to B \tn A$ for each objects of $\C$, satisfying the well-known hexagon properties. Remark that, the category of vector spaces is a braided monoidal category with a trivial braiding $ \mathcal{R}_{V,W} \colon v \tn w \in V \tn W \mapsto w \tn v  \in W \tn V$. 
%\end{Example}

\begin{remark}
Moreover, we have the following following properties and conventions:
\begin{itemize}
	\item A Hopf algebra $H$ is said to be cocommutative, iff $\D = \mathcal{R} \circ \D$.
	\item A Hopf algebra morphism is a bialgebra morphism compatible with the antipode. 
	\item We use Sweedler's notation \cite{SW1} for coproducts; namely $\Delta (x) = \underset{(x)}{\sum } \,\, x' \otimes x''$, for all $x \in H$.
	\item We use the notation $\nabla (x \tn y) = x y$ in the calculations, for the sake of simplicity.
	\item Let $\C$ be the category of vector spaces and $H$ be a Hopf algebra in $\C$. An element $x\in H$ is said to be ``primitive", if $\Delta(x)=x \tn 1 +1 \tn x$; and ``group like", if $\Delta(x)=x \tn x$. If $x$ is group-like then $\epsilon(x)=1$ and $S(x)=x^{-1}$; and if $x$ is primitive then $\epsilon(x)=0$ and $S(x)=-x$. The set of primitive elements $Prim(H)$ defines a Lie algebra, and the set of group-like elements $Gl(H)$ defines a group. Thus we have the functors:
	\begin{align*}
	\xymatrix@R=40pt@C=30pt{ \mathrm{\{Lie\ Algebras\}}  &
		\mathrm{\{Hopf\ Algebras\}} \ar[r]^{Gl} \ar[l]_{Prim} &
		\mathrm{\{Groups\}} \, .
	}
	\end{align*}
	From this point of view, Hopf algebras can be thought as a unification of groups and Lie algebras.
	
	\item 	Let $H$ be a Hopf algebra. A sub-Hopf algebra $A \subseteq H$ is a sub-object $A$, such that $\nabla( A \tn A) \subseteq A$, $\Delta(A) \subseteq A \tn A$ and $\eta(\kappa) \subseteq A$. 
	
	\item The category of Hopf algebras has a zero object $\kappa$, see \cite{AD1}. We have unique morphisms $\eta_A \colon \k \to A$ and $\e_A \colon A \to \k$. The zero morphism between two Hopf algebras $A$ and $B$ is therefore $\zeta_{A,B}=\z\colon A \to B$, where $z_{A,B}=\eta_B \circ \e_A$, thus $\z_{a,b}(x)=\e(x)1_B$.  
	%Clearly a sub-Hopf algebra inherits a Hopf algebra structure from $H$.
	%Moreover, a sub-Hopf algebra $A$ of $H$ is called normal \cite{vespa}, if $x \ad a \in A$ for all $x \in H$ and $a \in A$. 
\end{itemize}
\end{remark}

\subsection{Yetter-Drinfeld modules}

In this section, we recall the category of Yetter-Drinfeld modules which is defined via Hopf algebra actions and coactions. Moreover, this category is an essential example of a braided monoidal category.
   
\subsubsection{Modules and comodules}

\begin{Definition}[module]
Let $H$ be a bialgebra. $V$ is said to be a $H$-module (on the left) if we are given a bilinear map $\rho$ called (left) action: $$\rho\colon (x,v) \in H \times V \mapsto x \tr v \in V,$$ such that, for all $x,y \in H$ and each $v \in V$, we have:
$$(xy)\tr v=x\tr (y\tr v), \quad \quad 1_H \tr v=v,$$
where the first condition means the following diagram commutes:
\begin{align}\label{actioncond}
\xymatrix@R=20pt@C=20pt{ H \tn H \tn V  \ar[r]^-{\nabla \tn \id} \ar[d]_{\id \tn \rho} & H \tn V \ar[d]^{\rho} \\ H \tn V \ar[r]^-{\rho} & V }
\end{align} %We will denote by the same letter $\rho\colon H \times V \to V$, and also the associated bilinear map $\rho\colon H \tn V \to V$. 
\end{Definition}

\begin{Definition}[comodule]
	Let $H$ be a bialgebra. $V$ is said to be a (left) $H$-comodule, if we are given a linear map $\phi\colon V \to H \tn V$ called (left) coaction:
	$$\phi: v \in V \mapsto \sum_{[v]} v_H \tn v_V   \in H \tn V \, ,$$
	such that, for all $v \in V$, we have:
	\begin{align*}
	\sum_{(v_H)} \sum_{[v]} v_H' \tn v_H'' \tn v_V= \sum_{[v_V]} \sum_{[v]} v_H \tn (v_V)_H \tn  (v_V)_V \, ,
	\end{align*}
	namely the diagram commutes:
	\begin{align}\label{coactioncond}
	\xymatrix@R=20pt@C=20pt{ V  \ar[r]^-{\phi} \ar[d]_{\phi} & H \tn V \ar[d]^{\D \tn \id} \\ H \tn V \ar[r]^-{\id \tn \phi} & H \tn H \tn V } 
	\end{align}
	% A comodule is said to be invertible if we can define a linear map $k\colon V \to H \tn V$ such that:
	% $$\sum_{[v]} k(v_V) \left (v_H \tn v_V)=1 \tn V= $$
\end{Definition}

Analogously, one can define the right (co)action and therefore right (co)modules. However, we only use left (co)action and left $H$-(co)module conventions throughout this paper.

\subsubsection{Yetter-Drinfeld modules}

%Let $H$ be a Hopf algebra. Recall that $V$ is said to be a left $H$-comodule, if we are given a linear map $\phi\colon V \to H \tn V$ called left coaction:
%$$\phi: v \in V \mapsto \sum_{[v]} v_H \tn v_V   \in H \tn V \, ,$$
%such that, for all $v \in V$:
%\begin{align*}
%\sum_{(v_H)} \sum_{[v]} v_H' \tn v_H'' \tn v_V= \sum_{[v_V]} \sum_{[v]} v_H \tn (v_V)_H \tn  (v_V)_V \, ,
%\end{align*}
%namely the diagram commutes:
%\begin{align}\label{coactioncond}
%\xymatrix@R=20pt@C=20pt{ V  \ar[r]^-{\phi} \ar[d]_{\phi} & H \tn V \ar[d]^{\D \tn \id} \\ H \tn V \ar[r]^-{\id \tn \phi} & H \tn H \tn V } 
%\end{align}
% A comodule is said to be invertible if we can define a linear map $k\colon V \to H \tn V$ such that:
% $$\sum_{[v]} k(v_V) \left (v_H \tn v_V)=1 \tn V= $$

\begin{Definition}
Let $H$ be a bialgebra. We say that $V$ is a left-left Yetter-Drinfeld module\footnote{From now on, we simply call it a Yetter-Drinfeld module, since we only use left (co)actions. In the literature, Yetter-Drinfeld modules are also called ``crossed bimodules" and ``Yang-Baxter modules" \cite{zbMATH06812675}.} \cite{zbMATH00427768}, if:
\begin{itemize}
	\item $V$ is a left $H$-module,
	\item $V$ is a left $H$-comodule,
	\item The following compatibility condition holds:
	\begin{align}\label{comp2}
	\sum_{(x)} \sum_{[v]} x' v_H \tn x'' \t v_V= \sum_{(x)}\sum_{[x' \t v]} (x' \t v)_H \,\, x'' \tn (x' \t v)_V,
	\end{align}
	which is to say, the following diagram commutes:
	\begin{align}\label{compatibilitycond}
	\xymatrix@R=17pt@C=35pt{ H \tn V \ar[r]^-{\D \tn \id} \ar[d]_-{\D \tn \id} & H \tn H \tn V \ar[r]^-{\id \tn \id \tn \phi}  & H \tn H  \tn H \tn V \ar[dd]^-{\id \tn \mathcal{R} \tn \id}  \\
		H \tn H \tn V \ar[d]_-{\id \tn \mathcal{R}} & &  \\
		H \tn V \tn H \ar[d]_-{\rho \tn \id} &  & H \tn H \tn  H \tn V \ar[dd]^-{\nabla \tn \rho} \\
		V \tn H \ar[d]_-{\phi \tn \id} & &  \\
		H \tn V \tn H \ar[r]^-{\id \tn \mathcal{R}}  & H \tn H \tn V \ar[r]^-{\nabla \tn \id}  & H \tn V   } 
	\end{align}
\end{itemize}
\end{Definition}
%a left $H$-module and left $H$-comodule, such that, denoting the action of $H$ on $V$ by: $$\rho\colon x \tn v \in H \tn V \mapsto x \t v \in V,$$
% $$\phi: v \in V \mapsto \sum_{[v]} v_H \tn v_V   \in H \tn V$$
% therefore:
% $$\forall x , y \in H , \forall v \in V: (xy) \t v = x \t ( y \t v),  $$
% and
% $$\forall v \in V: \sum_{[v]} \sum_{(v_H)} v_H' \tn v_H'' \tn v_V= \sum_{[v]} \sum_{[v_V]} v_H \tn (v_V)_H \tn  (v_V)_V,$$
%the following compatibility condition holds:
%\begin{align}\label{comp2}
%\sum_{(x)} \sum_{[v]} x' v_H \tn x'' \t v_V= \sum_{(x)}\sum_{[x' \t v]} (x' \t v)_H \,\, x'' \tn (x' \t v)_V.
%\end{align}
%In other words, the following diagram commutes:
%\begin{align}\label{compatibilitycond}
%\xymatrix@R=25pt@C=50pt{ H \tn V \ar[r]^-{\D \tn \id} \ar[d]_-{\D \tn \id} & H \tn H \tn V \ar[r]^-{\id \tn \id \tn \phi}  & H \tn H  \tn H \tn V \ar[dd]^-{\id \tn \mathcal{R} \tn \id}  \\
%	H \tn H \tn V \ar[d]_-{\id \tn \mathcal{R}} & &  \\
%	H \tn V \tn H \ar[d]_-{\rho \tn \id} &  & H \tn H \tn  H \tn V \ar[dd]^-{\mu \tn \rho} \\
%	V \tn H \ar[d]_-{\phi \tn \id} & &  \\
%	H \tn V \tn H \ar[r]^-{\id \tn \mathcal{R}}  & H \tn H \tn V \ar[r]^-{\mu \tn \id}  & H \tn V   } 
%\end{align}

\begin{Notation}
	The category of Yetter-Drinfeld modules over a Hopf algebra $H$ in $\C$ will be denoted by $\mathcal{YD}_{\C}(H)$. Especially, if $\C$ is the category of vector spaces, we simply denote it by $\mathcal{YD}(H)$.
	%This notation will be specifically used when $H$ is a Hopf algebra in the category of vector spaces. 
\end{Notation}

\begin{Lemma}
Suppose that $H$ has an invertible antipode. Then, $\mathcal{YD}_{\C}(H)$ is braided. Here, given two Yetter-Drinfeld modules $V$ and $W$, we have the braiding:
\begin{align}\label{prebraiding}
\mathcal{R'} \colon v \tn w \in V \tn W \mapsto \sum_{[v]}v_H \t w \tn v_V \in W \tn V \, ,
\end{align}
which is obtained by the following diagram:
\begin{align}\label{prebraiding-dia}
\xymatrix@R=40pt@C=50pt{ V \tn W  \ar[r]^-{\phi \tn \id} \ar[d]_{\mathcal{R'}} & H \tn V \tn W \ar[d]^{\id \tn \mathcal{R}} \\ W \tn V \ar@{<-}[r]^-{\rho \tn \id} & H \tn W \tn  V } 
\end{align}
\end{Lemma}

\begin{remark}
If the antipode is not invertible, \eqref{prebraiding} is called a ``prebraiding" yielding a prebraided monoidal category in which the braiding map is not an isomorphism. For a detailed discussion about the alternative conditions to make $\mathcal{YD}_{\C}(H)$ into a braided monoidal category, see \cite{zbMATH00427768}. For this reason, all Hopf algebras will be considered to have an invertible antipode in the rest of the paper.
\end{remark}

\begin{Definition}[Adjoint Action]
	Let $A$ be a Hopf algebra. The left and right adjoint actions of $A$ on itself are:
	%\begin{align*}
	%\xymatrix@R=20pt@C=30pt{ A \tn A   \ar[r]^-{\mathrm{\D \tn \id}} \ar@/_1pc/[drrr]^{\ad} & A \tn A \tn A \ar[r]^-{\id \tn \mathcal{R}} & A \tn A \tn A  \ar[r]^-{\id \tn S \tn \id} & A \tn A \tn A \ar[d]^{\mathrm{\nabla \circ \nabla}} \\  &  & & A }  
	%\end{align*}
	\begin{align} \label{adjoint}
	a \ad b = \sum_{(a)} a' b S(a''), \quad \quad
	b \tl_{ad} a=\sum_{(a)} S(a') b a'' .
	\end{align}
\end{Definition}

\begin{Proposition}
	Suppose that $\d \colon I \to H$ is an algebra morphism. Then $H$ has a natural $I$-module algebra (Definition \ref{algebra-coalgebra}) structure with the action $\rho \colon I \tn H \to H$ given by:
	 \begin{align*}
	\xymatrix@R=30pt@C=15pt{
		I \tn H \ar[rr]^-{{\rho}} \ar[dr]_-{(\d \tn \id)}
		&& H \\ 	& H \tn H  \ar[ur]_-{\ad}
		& 
	} \end{align*}
Consequently, every Hopf algebra $A$ turns into an $A$-module algebra. On the other hand, note that:
	\begin{align*}
	\Delta(a \ad b) & = \sum_{(a)}\sum_{(b)} a' b' S(a'''')\tn a'' b'' S(a''') ,\\
	\Delta(b \tl_{ad} a) & = \sum_{(a)}\sum_{(b)} S(a'') b' a'''\tn S(a') b'' a'''' .
	\end{align*}
	From here we can see that neither adjoint action ought to turn, in general, $A$ into a module coalgebra over itself. This is however the case if $A$ is cocommutative.
\end{Proposition}

\begin{Definition}[Adjoint Coaction]
	The left and right adjoint coactions of a Hopf algebra on itself are:
	\begin{align}\label{coadjoint}
	\phi{(y)}=\sum_{(y)} y' S(y''') \tn y'', \quad \quad \phi{(y)}=\sum_{(y)} y'' \tn S(y') y''',
	\end{align}
	which are obtained in a dual way to adjoint actions as:
	\begin{align}\label{adj-coadj-dia}
	\xymatrix@R=30pt@C=40pt{ H \tn H \ar[d]_-{\D \tn \id} \ar@{-->}[rr]^-{\eqref{adjoint}} & &  H \ar@/^-2.5pc/[dddrr]_-{\D \circ \D}  \ar@{-->}[rr]^-{\eqref{coadjoint}} & &  H \tn H \\
		H \tn H \tn H \ar[d]_-{\id \tn \mathcal{R}} & &  & & H \tn H \tn H  \ar[u]_-{\nabla \tn \id}^-{\id \tn \nabla}  \\ 
		H \tn H \tn H \ar[d]_-{\id \tn \id \tn S}^-{S \tn \id \tn \id}  & & & & H \tn H \tn H \ar[u]_-{\id \tn \mathcal{R}} \\ 
		H \tn H \tn H  \ar@/^-2.5pc/[uuurr]_-{\nabla \circ \nabla} & & & & H \tn H \tn H \ar[u]_-{\id \tn \id \tn S}\\ } 
	\end{align}
\end{Definition}

\begin{Proposition}
	Suppose that $\d \colon I \to H$ is a coalgebra morphism. Then $H$ has a natural $I$-comodule coalgebra (Definition \ref{algebra-coalgebra}) structure with the coaction $\rho \colon I \to I \tn H$ given by:
	\begin{align*}
	\xymatrix@R=30pt@C=15pt{
		I \ar[rr]^-{{\rho}} \ar[dr]_-{(\D)}
		&& I \tn H \\ 	& I \tn I  \ar[ur]_-{\eqref{coadjoint}}
		& 
	} \end{align*}
	Consequently, every Hopf algebra $A$ turns into an $A$-comodule coalgebra.
\end{Proposition}

\begin{Lemma}
 Let $H$ be a Hopf algebra. Then $H$ itself is a Yetter-Drinfeld module with the adjoint action \eqref{adjoint} and regular coaction $\phi =\Delta$.
 \end{Lemma}
The previous construction has a natural extension by considering Hopf algebra projections as follows:

\begin{Definition}[Hopf algebra projection]
A Hopf algebra projection $(\d\colon I \to H,i)$ is given by a Hopf algebra morphism $\d$ together with another Hopf algebra morphism $i\colon H \to I$ (an inclusion) such that $\d\, i=\id_H$.
\end{Definition}
\begin{Lemma}\label{proj-lemma}
 Let $(\d\colon I \to H,i)$ be a Hopf algebra projection. Then, $I$ has a natural Yetter-Drinfeld module structure over $H$ with the action $\rho \colon H \tn I \to I$ and coaction $\phi \colon I \to H \tn I$ given by:
 \begin{align*}
 \xymatrix@R=30pt@C=15pt{
 	H \tn I \ar[rr]^-{{\rho}} \ar[dr]_-{(i \tn \id)}
 	&& I \\ 	& I \tn I  \ar[ur]_-{\ad}
 	& 
 } \quad \quad \quad \quad \quad
 \xymatrix@R=30pt@C=15pt{
 	I \ar[rr]^-{{\phi}} \ar[dr]_-{\D}
 	&& H \tn I \\ 	& I \tn I  \ar[ur]_-{(\d \tn \id)}
 	& 
 }
 \end{align*}
% \begin{align*}
% \rho \colon & H \tn I \to I \quad \text{by} \quad x \t v = i(x) \tad v \, , \\
% \phi \colon & I \to H \tn I \quad \text{by} \quad \phi = (\d \tn \id) \circ \Delta \, . 
% \end{align*}
\end{Lemma}

\subsubsection{Braided Hopf algebras in $\mathcal{YD}(H)$}\label{braided-in-yd}

The main idea of this part is to obtain a Hopf algebra structure living in the category of Yetter-Drinfeld modules. For this reason, we follow the categorical meaning of braided Hopf algebras \cite{zbMATH01824088,Rad1}. However, there is also a non-categorical definition of braided Hopf algebras given in \cite{zbMATH01584734}.

\begin{Definition}\label{algebra-coalgebra}
	Let $H$ be a Hopf algebra. A braided Hopf algebra living in $\mathcal{YD}_{\C}(H)$ is given by:
	\begin{enumerate}
		\item An object $A$ of $\mathcal{YD}_{\C}(H)$. 
		\item A unital associative algebra structure $(A,\nabla,\eta)$, making it 
		\begin{itemize}
			\item a $H$-module algebra with: $$x \t (ab)=\sum_{x} (x' \t a) \, (x'' \t b), \quad \an \quad x \t 1_A=\epsilon(x) 1_A. $$
			That means, $A$ is an $H$-module and also a $\k$-algebra, and $\nabla,\eta$ are module maps.
			\begin{align}
			\xymatrix@R=30pt@C=40pt{ H \tn A \tn A  \ar[r]^-{\id \tn \nabla} \ar[d]_{\D \tn \id} & H \tn A \ar[r]^-{\rho} & A \ar@{<-}[d]^{\nabla} \\ H \tn H \tn A \tn A  \ar[r]^-{\id \tn \mathcal{R} \tn \id} & H \tn A \tn H \tn A \ar[r]^-{\rho \tn \rho} & A \tn A } 
			\end{align}
			\item  a $H$-comodule algebra with: $$ \phi(ab)= \sum_{[ab]} (ab)_H \tn (ab)_A = \sum_{[a]}\sum_{[b]} a_H b_H \tn a_A b_A \quad \an \quad \phi(1)=1_H \tn 1_A. $$
			That means, $A$ is an $H$-comodule and also a $\k$-algebra, and $\nabla,\eta$ are comodule maps.
			\begin{align}
			\xymatrix@R=30pt@C=30pt{ A \tn A  \ar[r]^-{\nabla} \ar[d]_{\phi \tn \phi} & A \ar[r]^-{\rho} & H \tn A \ar@{<-}[d]^{\nabla \tn \nabla} \\ H \tn A \tn H \tn A  \ar[rr]^-{\id \tn \mathcal{R} \tn \id} &  & H \tn H \tn A \tn A } 
			\end{align}
		\end{itemize}
		\item A counital coassociative coalgebra structure $(A,\bD,\e)$ with the notation of coproduct:  $$\bD(a)=\sum_{|a|} \ba \tn \bba.$$ making it
		\begin{itemize}
			\item a $H$-module coalgebra with: 
			$$\dis\sum_{|x \t a|} \underline{x \t a} \tn  \underline{\underline{x \t a}}= \sum_{(x)}\sum_{|a|} x' \t \underline{a}
			\tn x'' \t \underline{\underline{a}} \quad \an \quad \e(x \t a)=\e(x) \e(a).$$
			That means, $A$ is an $H$-module and also a $\k$-coalgebra, and $\bD,\e$ are module maps.
			\begin{align}
			\xymatrix@R=30pt@C=30pt{ H \tn A  \ar[r]^-{\rho \tn \id} \ar[d]_{\D \tn \bD} & A \ar[r]^-{\bD} & A \tn A \ar@{<-}[d]^{\rho \tn \rho} \\ H \tn H \tn A \tn A  \ar[rr]^-{\id \tn \mathcal{R} \tn \id} & & H \tn A \tn H \tn A } 
			\end{align}
			
			\item  a $H$-comodule coalgebra with:
			$$\dis \sum_{[a]} \sum_{|a_A|} a_H \tn \underline{a_A} \tn \underline{\underline{a_A}}=\sum_{|a|}\sum_{[\underline{a}]}\sum_{[\underline{\underline{a}}]} \underline{a}_H \underline{\underline{a}}_H \tn \underline{a}_A \tn \underline{\underline{a}}_A \quad \an \quad  \sum_{[a]} a_H \e(a_A)=\e(a)1_A.$$
			That means, $A$ is an $H$-comodule and also a $\k$-coalgebra, and $\bD,\e$ are comodule maps.
			\begin{align}
			\xymatrix@R=30pt@C=30pt{ A  \ar[r]^-{\phi} \ar[d]_{\bD} & H \tn A \ar[r]^-{\id \tn \bD} & H \tn A \tn A \ar@{<-}[d]^{\nabla \tn \id} \\ A \tn A  \ar[r]^-{\phi \tn \phi} & H \tn A \tn H \tn A \ar[r]^-{\id \tn \mathcal{R} \tn \id} & H \tn H \tn A \tn A } 
			\end{align}
		\end{itemize}
	
		\item A unital algebra morphisms $\bD\colon A \to A \otimes_{H} A$ and $\epsilon\colon A \to \k$. Here the product in $A\tn_H A$ is:
		$$(a \tn b) (c \tn d)=\sum_{[b]} (a \, b_H \t c )\tn ( b_A\, d), $$
		namely:
		\begin{align*}
		\xymatrix@R=30pt@C=45pt{ A \tn A \tn A \tn A  \ar[rr]^-{} \ar[d]_{\id \tn \phi \tn \id \tn \id} & & A \tn A \ar@{<-}[d]^{\nabla \tn \nabla} \\ A \tn H \tn A \tn A \tn A  \ar[r]^-{\id \tn \id \tn \mathcal{R} \tn \id} & A \tn H \tn A \tn A \tn A  \ar[r]^-{\id \tn \rho \tn \id \tn \id} & H \tn A \tn H \tn A } 
		\end{align*}
		\item An antipode $\bS\colon A \to A$ satisfying the usual properties given in \eqref{antipode}.
%		$\mathrm{(H4)}$, namely:
%		$$\sum_{|a|} \bS(\ba) \bba=\e(a) 1_A= \sum_{|a|} \ba\bS(\bba) \quad \an \quad \sum_{|a|} \epsilon(\ba) \bba=a= \sum_{|a|} \ba\epsilon(\bba) $$
%		and moreover: %(consequence of the axioms; see Majid's paper).
%		$$\bS(ab)=\nabla\big(\mathcal{R'}(\bS(a)\tn \bS(b))) $$
%		$$\bS(1_A)=1_A $$
%		$$(\bS \tn \bS)\circ \mathcal{R'} \circ \bD = \bD \circ \bS. $$
	\end{enumerate}
\end{Definition}

\begin{Notation}
	A braided Hopf algebra $A$ living in the category of Yetter-Drinfeld modules over $H$ will be called a $\mathcal{YD}(H)$-Hopf algebra. 
\end{Notation}

\begin{Definition}
Let $A$ be a $\mathcal{YD}(G)$-Hopf algebra, and $B$ be a $\mathcal{YD}(H)$-Hopf algebra.
A (braided) map between $t \colon A \to B$ consists of a Hopf algebra morphism $r \colon G \to H$ such that $t$ is:
\begin{itemize}
\item an algebra morphism (to preserve products),
\item a coalgebra morphism (to preserve coproducts),
\item a Yetter-Drinfeld module morphism to preserve the actions and coactions; namely the diagrams below commute:
 $$\xymatrix{G \tn A  \ar[r]^{r \tn t} \ar[d]_{\rho_A} & H \tn B\ar[d]^{\rho_B} \\ A\ar[r]^t & B } \quad \quad \quad
 \xymatrix{G \tn A  \ar[r]^{r \tn t} \ar@{<-}[d]_{\phi_A} & H \tn B\ar@{<-}[d]^{\phi_B} \\ A\ar[r]^t & B } $$
\end{itemize}
\end{Definition}

\subsubsection{More on module algebras}

\n The following is well known and a proof appears in \cite{SM1}.
\begin{Theorem}[Majid]\label{majmain0}
	If $I$ is an $H$-module algebra with the action $\rho$, then we can define an algebra $I \otimes_\rho H$ with the underlying vector space $I \otimes H$, with product:
	$$(u \tn x)(v \tn y)=\sum_{(x)}\big( u \,\, x' \tr v \big)\tn\big( x'' y \big), \we u,v \in I \an x,y \in H, $$
	and identity $1_I \otimes 1_H$. If $H$ is an $H$-module coalgebra, additionally satisfying the following compatibility\footnote{The compatibility condition \eqref{comp} is used to prove that the coproduct of $I \otimes_\rho H$ is an algebra morphism.} condition:
	\begin{equation}\label{comp}
	\sum_{(x)} x' \tn ( x''\tr v)=\sum_{(x)} x'' \tn ( x'\tr v), \textrm{ for each } x \in H \textrm{ and } v \in I, 
	\end{equation}
	then we have a bialgebra structure living in $I \otimes_\rho H$ with coproduct:
	$$\D(u \tn x)= \sum_{(u),(x)} (u'\tn x') \tn (u'' \tn x'') , \we x \in H \an u \in I,$$
	and: 
	$$\epsilon (u \tn x)= \e(u) \tn \e(x), \we x \in H \an u \in I.$$
	If, in addition, $H$ and $I$ are Hopf algebras, then $I \otimes_\rho H$ is a Hopf algebra with antipode:
	$$S(u \tn x)=\big ( 1_I \tn S(x) \big)\, (S(v) \tn 1_H),  \we x \in H \an u \in I.  $$
	%	Moreover the inclusions of $H$ and $I$ on  $I \otimes_\rho H$, namely:
	%	\begin{align*}\allowdisplaybreaks
	%	& i_H \colon x \in H  \mapsto 1_I \tn x \in I \otimes_\rho H, \quad  \quad
	%	i_I \colon v \in I \mapsto v \tn 1_I \in I \otimes_\rho H ,
	%	\end{align*}
	%	are bialgebra morphisms (therefore Hopf algebra morphisms if $H$ and $I$ are Hopf algebras). 
\end{Theorem}

\section{A Functorial Approach to Radford's Theorem}

We fix a braided monoidal category with equalizers $\C$ in the rest of the paper.

\medskip

Let $\Omega \colon I \to H$ be a Hopf algebra morphism in $\C$. Put: 
\begin{align}\label{rker-definition}
\rk_{\C}(\Omega)=\{v \in I:\sum_{(v)} v' \tn \Omega(v'')=v \tn 1\} \, .
\end{align}
We know from \cite{AD1} that:
\begin{itemize}
	\item $\rk_{\C}(\Omega)$ only defines a subalgebra of $A$ in $\C$, with the identity. Moreover, it is a left coideal, i.e. $$\D(\r(f,g)) \subseteq A \tn \r(f,g).$$ 
	\item $\rk_{\C}(\Omega)$ is not the actual kernel of $\Omega$ in the category of Hopf algebras.
	\item $\rk_{\C}(\Omega)$ is not closed under the antipode, since $S\big( \rk_{\C}(\Omega) \big) \subseteq \lk_{\C}(\Omega)$, where we put:
	\begin{align*}
	\lk_{\C}(\Omega)=\{v \in I:\sum_{(v)} \Omega(v') \tn v'' = 1 \tn v\} \, .
	\end{align*}
	\item $\rk_{\C}(\Omega)$ is invariant under the (left) adjoint action.
	\item In fact, $\rk_{\C}(\Omega)$ is given by $\r_{\C}(f,\eta_H \circ \e_I)$ where we use:
	\begin{align*}
	\r(f,g) & = \Big\{ a \in A : \sum_{(a)} a' \tn f(a'')= \sum_{(a)} a' \tn g(a'') \Big\} \, .
	\end{align*}
	%which defines a kernel in the category of vector spaces, but not in the Hopf algebra level. Similary, $\lk_{\C}(\Omega)$ too.
	\item However, more general case of $\rk_{\C}(\Omega)$ is given by:
	%denoted by $\ck_{\C}(\Omega)$ where:
\begin{align*}
\ck_{\C}(\Omega) =	\left \{ a \in A  : \sum_{(a)} a' \tn f(a'') \tn a'''' = \sum_{(a)} a' \tn 1 \tn a'' \right \} \, ,
\end{align*}
that correctly defines the kernel in the category of Hopf algebras. Remark that, we have: $$\rk_{\C}(\Omega)=\ck_{\C}(\Omega)=\lk_{\C}(\Omega) \, ,$$ if the category is cocommutative.
\end{itemize}
For more details on the categorical properties of Hopf algebras, we refer \cite{Agore,AgoreII,SB2,zbMATH06005967}.

\subsection{Radford's Theorem}

Radford/Majid's theorem gives an idea to put a Hopf algebra structure on \eqref{rker-definition} when it is obtained from a Hopf algebra projection. The method was introduced in \cite{SM3,Rad1} in the category of vector spaces. Afterwards, it is generalized to any braided monoidal category in \cite{zbMATH01036413}.

%\bigskip
%
%Let $\Omega \colon I \to H$ be a Hopf algebra map in $\C$. We already know that
%\begin{align*}
%\rk_{\C}(\Omega)=\{v \in I:\sum_{(v)} v' \tn \Omega(v'')=v \tn 1\}
%\end{align*}
%is a subalgebra of $I$ and also a left-coideal, from Lemma \eqref{Lker}. 
%Moreover it is invariant under the adjoint action of $I$, from Lemma \eqref{rker-invariant}.
%Given the fact that:
%$$\D( u \tad v)=\sum_{(u)} \sum_{(v)} v' u' S(v'''') \tn v'' u'' S(v''')  ,$$
%we see that $\rk(\Omega)$ is invariant under the adjoint action of $I$. In other words, if $v \in I$ and $a \in \rk(\Omega)$ then $v \tad a \in \rk(\Omega)$.

\subsubsection{$\rk_{\C}(\d)$ as a braided Hopf algebra}

\begin{remark}
	Let $(\d\colon I \to H, i)$ be a Hopf algebra projection in $\C$. Then, $I$ has a natural Yetter-Drinfeld module structure over $H$, from Lemma \ref{proj-lemma}. Since $\rk_{\C}(\Omega)$ is invariant under the (left) adjoint action of $I$, it follows that $\rk_{\C}(\Omega)$ is, itself, a Yetter-Drinfeld module over $H$.
\end{remark}

\begin{Notation}
	We fix $\B = \rk_{\C}(\d)$ for a given Hopf algebra morphism $\d\colon I \to H$ in the rest of this section.
\end{Notation}

\begin{Definition}
	For a given Hopf algebra projection $(\d\colon I \to H, i)$ we define maps\footnote{These maps are not Hopf algebra morphisms. Moreover, we call them ``kernel generator" maps.} $f,g \colon I \to I$ as follows:
	 \begin{align}\label{kernel-generators}
	\xymatrix@R=40pt@C=40pt{ I  \ar@{-->}[r]^-{f} \ar[d]_{\D} & I  \\ I \tn I \ar[r]^-{\id \, \tn \, i \d S} & I \tn I \ar[u]_{\nabla}} \quad \quad \quad
	\xymatrix@R=40pt@C=40pt{ I  \ar@{-->}[r]^-{g} \ar[d]_{\D} & I  \\ I \tn I \ar[r]^-{i \d \, \tn \, S} & I \tn I \ar[u]_{\nabla}}
	\end{align}
%	$$f\colon v \in  I \mapsto \sum_{(v)} v'\, i(\d(S(v'')) \in I, $$
%	$$g\colon v \in  I \mapsto \sum_{(v)} i(\d(v'))\, S(v'') \in I. $$
\end{Definition}

\begin{remark}
	Both $f,g$ take their values in $\B = \rk_{\C}(\d)$. Moreover, $f_{|\B}=\id_\B$, therefore $f^{2}=f$.
	% $I \to A$.
	On the other hand, we have: $$f(v\,i(x))=f(v) \an f \big( i(x) \ad v \big)= f \big( i(x) \, v \big)=x \ad f(v),$$
and also $g f=g $. Therefore:
\begin{align*}
\sum_{(v)} f(v')g(v'')=\epsilon(v), \quad \quad \sum_{(v)} g\big( f(v')\,v'' \big)=\epsilon(v), \quad \quad \sum_{(v)} f(v')i(\d(v''))=v.
\end{align*}
\end{remark}

\begin{Lemma}
	Let $(\d\colon I \to H,i)$ be a Hopf algebra projection in $\C$. Following the previous constructions, $I$ forms a $\mathcal{YD}_{\C}(H)$-Hopf algebra structure with the following new coproduct and antipode:
	\begin{align}\label{structure}
	\bD=(f \tn \id) \circ \D \, , \quad \quad \bS=g. 
	\end{align}
	
	In general, we know that $\B=\rk_{\C}(\d)$ is not a sub-Hopf algebra of $I$ in $\C$. However, since both $f$ and $g$ take values in $\B$, and $\B$ is invariant under the adjoint action, $\B$ forms a  sub-$\mathcal{YD}_{\C}(H)$-Hopf algebra of $I$. Consequently, $\B$ is also a $\mathcal{YD}_{\C}(H)$-Hopf algebra.
\end{Lemma}

\begin{Definition}[Braided adjoint action]
	Suppose that $A$ is a $\mathcal{YD}_{\C}(H)$-Hopf algebra. The braided adjoint action of $A$ on itself is given by:
	\begin{align}\label{bad}
	a \bad b= \sum_{|a|}\sum_{[\bba]} \ba \,\, (\bba_H \t b )\,\, \bS( \bba_A) \, .
	\end{align}
	Notice that, this is obtained in a same diagrammatic way to normal adjoint action given in \eqref{adj-coadj-dia}. However, the braiding map is not trivial in the category $\mathcal{YD}(H)$ which makes the difference in the formulae. 
	%In the case when $A$ arises from a Hopf algebra projection $(\d \colon I \to H,i)$, the braided adjoint action is just the adjoint action restricted to $A$.
\end{Definition}

\subsubsection{Bosonisation}

\begin{Lemma}\label{boso}
	Let $A$ be any $\mathcal{YD}_{\C}(H)$-Hopf algebra. We have a Hopf algebra $A \btn H$ with the underlying tensor product $A\tn H$, and with:
	$$(a \tn x) (b \tn y)=\sum_{(x)}(a \, x' \t b) \tn ( x'' y) \, ,$$
	$$\bD(a \tn x)=\sum_{|a|} \sum_{ [\bba]} \sum_{(x)} \ba \tn \bba_H x' \tn \bba_A \tn x'' ,$$
	$$S(a \tn x)=(1 \tn S(x))\,  (\bS(a) \tn 1)  \, .$$
	
	\medskip
	
	The smash product and the smash coproduct given above is categorically defined in a dual way to each other diagrammatically by:
	%\begin{align*}
	%	\xymatrix@R=40pt@C=40pt{ A \tn H \tn A \tn H \ar[d]^-{\id \tn \D \tn \id \tn \id} \ar@{-->}[rr]^-{\text{smash product}} & &  A \tn H \ar@/^-2.5pc/[dddrr]_-{\D \tn \D}  \ar@{-->}[rr]^-{\text{smash coproduct}} & &  A \tn H \tn A \tn H \\
	%	A \tn H \tn H \tn A \tn H \ar[d]^-{\id \tn \id \tn \tau \tn \id} & &  & & A \tn H \tn H \tn A \tn H  \ar[u]^-{\id \tn \mu \tn \id \tn \id}  \\ 
	%		A \tn H \tn A \tn H \tn H \ar[d]^-{\id \rho \tn \id \tn \id}  & & & & A \tn H \tn A \tn H \tn H \ar[u]^-{\id \tn \id \tn \tau \tn \id} \\ 
	%		A \tn A \tn H \tn H  \ar@/^-2.5pc/[uuurr]_-{\mu \tn \mu} & & & & A \tn A \tn H \tn H \ar[u]^-{\id \tn \phi \tn \id \tn \id}\\ } 
	%\end{align*}
	\begin{align*}
	\xymatrix@R=30pt@C=40pt{ 
		A \tn H \ar@{-->}[rr]^-{\text{smash coproduct}} \ar@/^-2.5pc/[dddrr]_-{\D \tn \D}  & & A \tn H \tn A \tn H \ar@<0.5ex>[d]^-{\id \tn \D \tn \id \tn \id} \ar@{-->}[rr]^-{\text{smash product}} & &  A \tn H  \\
		&& A \tn H \tn H \tn A \tn H \ar@<0.5ex>[d]^-{\id \tn \id \tn \tau \tn \id} \ar@<0.5ex>[u]^-{\id \tn \nabla \tn \id \tn \id} & &   \\ 
		&& A \tn H \tn A \tn H \tn H \ar@<0.5ex>[d]^-{\id \tn \rho \tn \id \tn \id} \ar@<0.5ex>[u]^-{\id \tn \id \tn \tau \tn \id}  & & \\ 
		&& A \tn A \tn H \tn H  \ar@/^-2.5pc/[uuurr]_-{\nabla \tn \nabla} \ar@<0.5ex>[u]^-{\id \tn \phi \tn \id \tn \id} & & } 
	\end{align*}
\end{Lemma}

\medskip

\begin{Theorem}
	Let $(\d\colon I \to H,i)$ be a Hopf algebra projection in $\C$. Then the pair of morphisms:
	$$\Psi\colon v \in I \mapsto \sum_{(v)} f(v') \tn \d(v'') \in \B \btn H.  $$
	$$\Phi\colon a \tn x \in \B \btn H \mapsto  a \, i(x) \in I$$
	are mutually inverse, defining the isomorphism of Hopf algebras $I \cong \B \btn H$ in $\C$.
\end{Theorem}

\begin{Corollary}
	Therefore, any $\mathcal{YD}_{\C}(H)$-Hopf algebra can arise from a Hopf algebra projection.
\end{Corollary}

\subsection{Interchanging Yetter-Drinfeld modules}\label{interchange}
 
\begin{Theorem}
Suppose that we have a Hopf algebra projection $ (\d \colon I \to H , i) $ in $\C$, and $B$ is any Yetter-Drinfeld Module over $H$. That means, we already have an action and coaction of $H$ on $B$ given by:
\begin{align*}
 \rho \, \colon & H \tn B \to B \, , \\
 \phi \, \colon & B \to H \tn B \, .
\end{align*}

Now, let us define $\bar{\rho} \colon I \tn B \to B$ and $\bar{\phi} \, \colon B \to I \tn B$ as follows:
\begin{align}\label{extended-actions}
\xymatrix@R=30pt@C=15pt{
	I \tn B \ar[rr]^-{\bar{\rho}} \ar[dr]_-{(\d \tn id)}
	&& B \\ 	& H \tn B  \ar[ur]_-{\rho}
	& 
} \quad \quad \quad \quad
\xymatrix@R=30pt@C=15pt{
	B \ar[rr]^-{\bar{\phi}} \ar[dr]_-{\phi}
	&& I \tn B \\ 	& H \tn B  \ar[ur]_-{(i \tn id)}
	& 
}
\end{align}
%\begin{align*}
% \bar{\rho} \, \colon & \bar{H} \tn B \to B \quad = \quad \bar{H} \tn B \xrightarrow{(\d \tn id)} H \tn B \xrightarrow{\,\,\,\, \rho \,\,\,\,} B \\
% \bar{\chi} \, \colon & B \to \bar{H} \tn B \quad = \quad B \xrightarrow{\,\,\,\, \chi \,\,\,\,} H \tn B \xrightarrow{(\i \tn id)} \bar{H} \tn B
%\end{align*}
With this (induced) action and coaction, $B$ also defines a Yetter-Drinfeld module over $I$. %\missing{Furthermore, it has (unital) module algebra structure in $\mathcal{YD}_{\C}(I)$.}
\end{Theorem}
 
 \medskip
 
 \begin{Proof}
 	First of all, it is clear that $1_I \tr b = b$ and $\phi(1)=1_I \tn 1_B$ for all $b \in B$ by definition, since $\d$ and $i$ are (Hopf) algebra morphisms. Moreover:
 	\begin{itemize}
 		\item $B$ is an $I$-module: In the diagram below, top square commutes since $\d$ is a Hopf algebra map. Also, it is easy to see that the bottom-left square also commutes. Moreover, since $B$ is an $H$-module, bottom-right square commutes, regarding \eqref{actioncond}. Consequently, the outer diagram commutes that proves $B$ is an $I$-module.
 		\begin{align*}
 		\xymatrix@R=40pt@C=50pt{ 
 			I \tn I \tn B \ar[rr]^-{(\nabla \tn \id)} \ar[d]_{(\id \tn \d \tn \id)} & & I \tn B \ar[d]^{(\d \tn \id)} \\ 
 			I \tn H \tn B \ar@{.>}[r]^-{(\d \tn \id \tn \id)} \ar[d]_-{(\id \tn \rho)} & H \tn H \tn B \ar@{.>}[r]^-{(\nabla \tn \id)} \ar@{.>}[d]^-{(\id \tn \rho)} & H \tn B \ar[d]^{\rho}  \\ 
 			I \tn B \ar[r]_-{(\d \tn \id)}  & H \tn B \ar[r]_-{\rho}  & B } 
 		\end{align*}

 		\item $B$ is an $I$-comodule: In the diagram below, right square commutes since $i$ is a Hopf algebra map, and also the bottom-left square already commutes. Since $B$ is an $H$-comodule, top-left square commutes, regarding \eqref{coactioncond}. Consequently, the outer diagram commutes that proves $B$ is an $I$-comodule.
 		\begin{align*}
 		\xymatrix@R=40pt@C=50pt{ 
 			B \ar[r]^-{\phi} \ar[d]_{\phi} & H \tn B \ar[r]^-{(i \tn \id)} \ar@{.>}[d]^-{(\D \tn \id)} & I \tn B \ar[dd]^{\D \tn \id} \\ 
 			H \tn B \ar@{.>}[r]^-{(\id \tn \phi)} \ar[d]_-{(i \tn \id)} & H \tn H \tn B  \ar@{.>}[d]^-{(i \tn \id \tn \id)} &  \\ 
 			I \tn B \ar[r]_-{(\id \tn \phi)}  & I \tn H \tn B \ar[r]_-{(\id \tn i \tn \id)}  & I \tn I \tn B } 
 		\end{align*}
 				
 		\item We check the compatibility condition: In the diagram below, the inner square commutes since $B$ is a YD-module over $H$, see \eqref{compatibilitycond}. It is readily checked that all other small diagrams commute. Therefore, the outer diagram commutes that proves the required compatibility.
 		\end{itemize}
  		\begin{align*}
 \xymatrix@R=50pt@C=35pt{  
 	I \tn B \ar[d]_-{(\D \tn \id)} \ar[rr]^-{(\D \tn \id)} & & I \tn I \tn B \ar[r]^-{(\id \tn \id \tn \phi)} & I \tn I \tn H \tn B \ar[r]^-{(\id \tn \id \tn i \tn \id)} & I \tn I \tn I \tn B  \ar[dd]^-{(\id \tn \mathcal{R} \tn \id)} \\ 
 	I \tn I \tn B  \ar[d]_-{(\id \tn \mathcal{R})}  & H \tn B  \ar[d]^-{(\D \tn \id)}   \ar[r]^-{(\D \tn \id)} \ar@{.>}[ul]_-{(i \tn \id)} & H \tn H \tn B \ar@{.>}[u]^-{(i \tn i \tn \id)} \ar[r]^-{(\id \tn \id \tn \phi)} & H \tn H \tn H \tn B \ar@{.>}[ur]^-{(i \tn i \tn i \tn \id)}  \ar[dd]_-{(\id \tn \mathcal{R} \tn \id)}  & \\
 	I \tn B \tn I  \ar[d]_-{(\d \tn \id \tn \id)}  & H \tn H \tn B  \ar[d]^-{(\id \tn \mathcal{R})}  \ar@{.>}[ul]_-{(i \tn i \tn \id)}  &&& I \tn I \tn I \tn B  \ar[dd]^-{(\nabla \tn \d \tn \id)}  \\
 	H \tn B \tn I  \ar[d]_-{(\rho \tn \id)} & H \tn B \tn H  \ar[d]^-{(\rho \tn \id)}   \ar@{.>}[l]_-{(\id \tn \id \tn i)}  & & H \tn H \tn H \tn B  
 	\ar@{.>}[ur]^-{(i \tn i \tn i \tn \id)}  \ar[dd]_-{(\nabla \tn \rho)} & \\
 	B \tn I  \ar[d]_-{(\phi \tn \id)} & B \tn H  \ar[d]^-{(\phi \tn \id)}  \ar@{.>}[l]_-{(\id \tn i)}  & & & I \tn H \tn B  \ar[dd]^-{(\id \tn \rho)}  \\
 	H \tn B \tn I  \ar[d]_-{(i \tn \id \tn \id)} & H \tn B \tn H  \ar@{.>}[dl]^-{(i \tn \id \tn i)}  \ar[r]^-{(\id \tn \mathcal{R})}  & H\tn H \tn B  \ar[r]^-{(\nabla \tn \id)} & H \tn B  \ar@{.>}[dr]_-{(i \tn \id)}  & \\
 	I \tn B \tn I  \ar[rr]_-{(\id \tn \mathcal{R})} & & I \tn I \tn B  \ar[rr]_-{(\nabla \tn \id)}  & & I \tn B
 } 
 \end{align*}
 \end{Proof}

\bigskip

\n Consequently:
\begin{Corollary}
	For a given Hopf algebra projection $ (\d \colon I \to H , i) $ in $\C$ and any object $B$ of $\yd$, we proved that $B$ can be converted to an object of $\mathcal{YD}_{\C}(I)$. Therefore, we obtain a monoidal functor:
	\begin{align*}
	\mathcal{F} \colon \yd \longrightarrow \mathcal{YD}_{\C}(I). 
	\end{align*}
%	Moreover, we can dualize the below functor in a similar way. In other words, if we have an object $A$ of $\mathcal{YD}_{\C}(I) $ in a Hopf algebra projection $ (\d \colon I \to H , i) $, $A$ also defines an object of $\yd$ where the action and coaction of $H$ on $A$ is given dual way to \eqref{extended-actions}. Consequently, we get two monoidal functors:
%	\begin{align}\label{functors}
%	\xymatrix{ \mathcal{YD}_{\C}(H)  \ar@<1ex>[rr]^{\mathcal{F}} && \mathcal{YD}_{\C}(I)  \ar@<1ex>[ll]^{\mathcal{G}}},
%	\end{align}
%	where $\mathcal{G} \circ \mathcal{F} = \id_{\mathcal{YD}_{\C}(H)}$ and also $\mathcal{F} \circ \mathcal{G} = id_{\mathcal{YD}_{\C}(I)}$.
\end{Corollary}

\begin{Proposition}
 The functor $\F$ is braided monoidal, since it is identical on the braiding.
\end{Proposition}
 
 \begin{Proof}
  Recalling the braiding $\mathcal{R'}$ in $\mathcal{YD}_{\C}(H)$ from \eqref{prebraiding-dia}, we write $\F (\mathcal{R'})$ as follows:
%  \begin{center}
%  	$\F \, \left( 
%  	\begin{tabular}{c}
%  	$
%  	\xymatrix@R=40pt@C=50pt{ V \tn W  \ar[r]^-{\phi \tn \id} \ar[d]_{\mathcal{R'}} & H \tn V \tn W \ar[d]^{\id \tn \mathcal{R_{\C}}} \\ W \tn V \ar@{<-}[r]^-{\rho \tn \id} & H \tn W \tn  V }
%  	$
%  	\end{tabular}
%  	\right) \quad = 
%  	\begin{tabular}{c}
%  	$\xymatrix@R=40pt@C=50pt{ \F (V \tn W)  \ar[r]^-{\F(\phi \tn \id)} \ar[d]_{\F(\mathcal{R'})} & \F (H \tn V \tn W) \ar[d]^{\F(\id \tn \mathcal{R_{\C}})} \\ \F (W \tn V) \ar@{<-}[r]^-{\F (\rho \tn \id)} & \F (H \tn W \tn  V) } $
%  	\end{tabular}$
%  \end{center}
%  which is equal to
   \begin{align}
  \xymatrix@R=40pt@C=50pt{ V \tn W  \ar[r]^-{\bar{\phi} \tn \id} \ar[d]_{\F (\mathcal{R'})} & I \tn V \tn W \ar[d]^{\id \tn \mathcal{R_{\C}}} \\ W \tn V \ar@{<-}[r]^-{\bar{\rho} \tn \id} & I \tn W \tn  V } 
  \end{align}
  which is (by definition):
    \begin{align*}
  \xymatrix@R=40pt@C=50pt{ V \tn W  \ar[r]^-{\phi \tn \id} \ar[d]_{\F (\mathcal{R'})} & H \tn V \tn W \ar[r]^-{i \tn \id \tn \id} & I \tn V \tn W \ar[d]^{\id \tn \mathcal{R_{\C}}}  \\ W \tn V \ar@{<-}[r]^-{\rho \tn \id} & H \tn W \tn  V & I \tn W \tn V \ar[l]_{\d \tn \id \tn \id} } 
  \end{align*}
  and using the fact that $\d \, i = \id$, we get:
    \begin{align}
  \xymatrix@R=40pt@C=50pt{ V \tn W  \ar[r]^-{\phi \tn \id} \ar[d]_{\F (\mathcal{R'})=\mathcal{R'}} & H \tn V \tn W \ar[d]^{\id \tn \mathcal{R_{\C}}} \\ W \tn V \ar@{<-}[r]^-{\rho \tn \id} & H \tn W \tn  V } 
  \end{align}
%  If we consider $\F (\mathcal{R'})$, we obtain:
%  \begin{align*}
%  \xymatrix@R=40pt@C=50pt{ \F (V \tn W)  \ar[r]^-{\F(\phi \tn \id)} \ar[d]_{\F(\mathcal{R'})} & \F (H \tn V \tn W) \ar[d]^{\F(\id \tn \mathcal{R_{\C}})} \\ \F (W \tn V) \ar@{<-}[r]^-{\F (\rho \tn \id)} & \F (H \tn W \tn  V) } 
%  \end{align*}
%  or, equivalently by definition \eqref{extended-actions}:
%  \begin{align*}
%  \xymatrix@R=40pt@C=50pt{ V \tn W  \ar[r]^-{\phi \tn \id} \ar[d]_{\F(\mathcal{R'})} & H \tn V \tn W \ar[r]^-{i \tn \id \tn \id} & I \tn V \tn W \ar[d]^{\id \tn \mathcal{R_{\C}}}  \\ W \tn V \ar@{<-}[r]^-{\rho \tn \id} & H \tn W \tn  V & I \tn W \tn V \ar[l]_{\d \tn \id \tn \id} } 
%  \end{align*}
%  which is equal to \eqref{prebraiding2} (remark that we used the fact that $\d \circ i = \id$). 
 \end{Proof}

%\begin{remark}
% Contrarily to theorem above, the functor $\mathcal{G}$ is not braided monoidal.
%\end{remark}

%\noindent Therefore it is clear that if $B$ is an object of category $\mathcal{DY}_{\C}(H)$, then $\F(B)$ is an object of category $\mathcal{DY}_{\C}(\bar{H})$. But we already know that, moreover, $B$ has a Hopf algebra structure in $\mathcal{DY}_{\C}(H)$. So:

\begin{Proposition}\label{main4}
	Suppose that we have a braided Hopf algebra structure $A$ living in  $\mathcal{YD}_{\C}(H)$. Then we have a natural Hopf algebra structure $\mathcal{F}(A)$ living in $ \mathcal{YD}_{\C}(I)$. 
%	with: 
%	\begin{align*}
%	\mathcal{H} (A) := \big( \mathcal{H}(A) , \mathcal{H}(\nabla) , \mathcal{H}(\eta) , \mathcal{H}(\D), \mathcal{H}(\epsilon), \mathcal{H}(S) \big) \, .
%	\end{align*}
%	where $\mathcal{H}(\nabla) \colon \mathcal{H}(A) \tn \mathcal{H}(A) \to \mathcal{H}(A)$ and $\mathcal{H}(\D) \colon \mathcal{H}(A) \to \mathcal{H}(A) \tn \mathcal{H}(A)$. 
\end{Proposition}

\begin{Proof}
	In fact, the functor $\mathcal{F}$ only changes the module structure, i.e. actions and coactions. Therefore, it does not change the object $A$ and Hopf algebra operations on it. The most crucial point here is, the braiding is used in the compatibility law $\eqref{compa}$. However, we also proved that the braiding is preserved (identical) under $\mathcal{H}$. 
\end{Proof}

\bigskip

\n Consequently, we can give the following:

\begin{Corollary}\label{main5}
	Suppose that we have:
	\begin{itemize}
		\item A Hopf algebra projection $ (\d \colon I \to H , i) $ in a braided monoidal category $\C$,
		\item A Hopf algebra $B$ living in $\yd$.
	\end{itemize}
	Then we obtain a Hopf algebra structure $\F(B)$ living in $\mathcal{YD}_{\C}(I)$.
\end{Corollary}

\section{Radford's Theorem in a Simplicial Structure}

In this section, we apply Radford's theorem to a simplicial Hopf algebra given in a braided monoidal category $\C$, based on the previous section.
%But first, let us recall the notion of simplicial object in any category.
%The cocommutative case of this problem is examined in \cite{1905.09620} for 2-crossed module of (cocommutative) Hopf algebras.

\subsection{Simplicial Hopf algebras}

A simplicial Hopf algebra $\mathcal{H}$ is a simplicial set in the category of Hopf algebras. In other words, it is a collection of Hopf algebras $%
H_{n}$ $(n\in \mathbb{N})$ in $\C$ together with Hopf algebra morphisms:
\begin{align*}
\begin{tabular}{llll}
$d_{i}^{n} \colon$ & $H_{n}\to H_{n-1}$ & $,$ & $0\leq i\leq n $ \\
$s_{j}^{n+1} \colon $ & $H_{n}\to H_{n+1}$ & $,$ & $0\leq j\leq n$%
\end{tabular}
\end{align*}
which are called faces and degeneracies, respectively\footnote{To avoid overloaded notation, we will not use superscripts for faces and degeneracies.}. These morphisms are requied to satisfy the following axioms, called ``simplicial identites'':%
\begin{equation*}
\begin{tabular}{rlll}
(1) & $d_{i}d_{j}=d_{j-1}d_{i}$ & if & $i<j$ \\
(2) & $s_{i}s_{j}=s_{j+1}s_{i}$ & if & $i\leq j$ \\
(3a) & $d_{i}s_{j}=s_{j-1}d_{i}$ & if & $i<j$ \\
(3b) & $d_{j}s_{j}=d_{j+1}s_{j}= \id$ &  &  \\
(3c)& $d_{i}s_{j}=s_{j}d_{i-1}$ & if & $i>j+1$ \\
\end{tabular}%
\end{equation*}

A simplicial Hopf algebra can be pictured as:
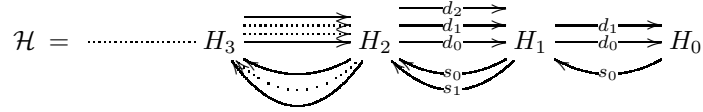
\begin{align*}
\xymatrix@C=40pt{
	\mathcal{H} \: = \: \: \ar@{.}[r] &
	{H}_{3}
	\ar@<2.25ex>[r] \ar@{.>}@<1.5ex>[r] \ar@{.>}@<0.75ex>[r] \ar@<0ex>[r] &
	{H}_{2}\ar@<3ex>[r]|{d_2} \ar@<1.5ex>[r]|{d_1} \ar@<0ex>[r]|{d_0}
	\ar@/^1pc/[l] \ar@{.>}@/^1.5pc/[l] \ar@/^2pc/[l]&
	{H}_{1}\ar@<1.5ex>[r]|{d_1} \ar[r]|{d_0}
	\ar@/^1pc/[l]|{s_0} \ar@/^1.5pc/[l]|{s_1} &
	{H}_{0}\ar@/^1pc/[l]|{s_0} }
\end{align*}

%\begin{Definition}[Truncation]
% An n-truncated simplicial Hopf algebra:
% \begin{align*}
%  Tr_{n}\mathcal{H}=\{H_{n},\ldots,H_{1},H_{0}\}
% \end{align*}
%  is a special simplicial Hopf algebra obtained by forgetting any dimensions of greater than n in the simplicial Hopf algebra $\mathcal{H}$. We denote category on n-truncated simplicial Hopf algebras by $Tr_{n}SimpHopf$, the full subcategory of $SimpHopf$.
%  
%  \medskip
%
%  Clearly, any simplicial Hopf algebra $\mathcal{H}$ yields an n-truncated simplicial Hopf algebra $Tr_{n}(\mathcal{H})$. Therefore we get the functor:
%  \begin{align*}
%   Tr_{n} \colon SimpHopf \to Tr_{n}SimpHopf
%  \end{align*}
%\end{Definition}

%\subsection{Radford's Theorem in a Simplicial Hopf Algebra}
%
%In this part, we apply Radford's Theorem in a simplicial structure step by step. This process gives rise to a proper definition of the Moore complex of a given simplicial Hopf algebra, which is one of the main outcomes of this paper. 

%Let $\mathcal{H}$ be a simplicial Hopf algebra given in \eqref{simplicial}
%$$\xymatrix@C=40pt{
%      \mathcal{H} \: = \: \: \ar@{.}[r] &
%      {H}_{3}
%      \ar@<2.25ex>[r] \ar@{.>}@<1.5ex>[r] \ar@{.>}@<0.75ex>[r] \ar@<0ex>[r] &
%      {H}_{2}\ar@<3ex>[r]|{d_2} \ar@<1.5ex>[r]|{d_1} \ar@<0ex>[r]|{d_0}
%      \ar@/^1pc/[l] \ar@{.>}@/^1.5pc/[l] \ar@/^2pc/[l]&
%      {H}_{1}\ar@<1.5ex>[r]|{d_1} \ar[r]|{d_0}
%      \ar@/^1pc/[l]|{i_0} \ar@/^1.5pc/[l]|{i_1} &
%      {H}_{0}\ar@/^1pc/[l]|{i_0} }$$

\begin{remark}
	In a simplicial Hopf algebra structure, we can obtain various Hopf algebra projections using the simplicial identity (3b), such as:
	\begin{align*}
	\big( d_{1} \colon H_{1} \to H_{0} \, , \, s_{0} \big) \, , \quad \big( d_{2} \colon H_{2} \to H_{1} \, , \, s_{1} \big) \, , \quad \big( d_{2} \colon H_{3} \to H_{2} \, , \, s_{2} \big) \, , \quad \text{etc.}
	\end{align*}
	Obviously, there exist $2n$ different Hopf algebra projections coming from (3b) between $H_n$ and $H_{n-1}$.
\end{remark}

\medskip

\begin{Notation}
	Suppose that we have a Hopf algebra projection:
	\begin{align*}
	\Big( \d_{j} \colon H_{n} \to H_{n-1} \, , \, s_{k} \Big) \, ,
	\end{align*}
	in a braided monoidal category $\mathcal{D}$. We put:
	\begin{align}\label{rker-gen}
	A_{j,k}^{n} = RKer_{\mathcal{D}} \big( \d_{j} \colon H_{n} \to H_{n-1} \big),
	\end{align}
	where the kernel generators \eqref{kernel-generators} will be denoted by:
	\begin{align}\label{aa}
	f^{n}_{j,k} = \nabla \big( \id \tn s_k d_j S \big) \D \, , \quad \quad g^{n}_{j,k}= \nabla \big(s_k d_j \tn S \big) \D \, ,
	\end{align}
	and following \eqref{rker-gen}, the coproduct of $A^{n}_{j,k}$ becomes:
	\begin{align}\label{bb}
	\D_{A^{n}_{j,k}} = \big(  f^{n}_{j,k} \tn \id \big) \D_{\C} \, ,
	\end{align}
	% \medskip
	% 
	% \centerline{\relabelbox
	% 	\epsfysize 2cm
	% 	\epsfbox{figs/coproduct-proj.eps}
	% 	\relabel{c}{$\scalebox{0.9}{$f^{\star}$}$}
	% 	\relabel{f}{$\scalebox{0.9}{$\D_{A^{n}_{j,k}}$}$}
	% 	\relabel{e}{$\scalebox{0.9}{$\D_{\C}$}$}
	% 	\relabel{b}{$\scalebox{0.8}{$=$}$}
	% 	\endrelabelbox}
	% 
	% \medskip
	% 
	and the antipode of $A^{n}_{j,k}$ is given by $g^{n}_{j,k}$ as usual.
\end{Notation}

\begin{Example}
	In a simplicial Hopf algebra, we have:
	\begin{align*}
	A_{0,0}^{2} = RKer_{\C} \big( d_{0} \colon H_{2} \to H_{1} \big) \, ,
	\end{align*}
	which is obtained from the Hopf algebra projection $\big( d_{0} \colon H_{2} \to H_{1} \, , \, s_{0} \big)$.
\end{Example}

%\subsection{Dimension one}
%
%Consequently, given a simplicial Hopf algebra \eqref{simplicial}, we have:
%\begin{align*}
%A_{0,0}^{n} = RKer_{\C} \big( d_{0} \colon H_{n} \to H_{n-1} \big) \, ,
%\end{align*}
%which are obtained from the Hopf algebra projection $\big( d_{0} \colon H_{n} \to H_{n-1} \, , \, i_{0} \big)$.

\subsection{Applying to dimension one and dimension two}\label{dimension2}

Consider the second part of the simplicial Hopf algebra, namely:
\begin{align}\label{h2-case}
\xymatrix@C=60pt@R=30pt{
	{H}_{2}\ar@<3ex>[r]|{d_2} \ar@<1.5ex>[r]|{d_1} \ar@<0ex>[r]|{d_0}
	& {H}_{1} \ar@<1.5ex>[l]|{s_0} \ar@<3ex>[l]|{s_1}  \\
	%	A_{0}^{2} \ar@<0.75ex>[r]^{d_1} \ar@<-0.75ex>@{<-}[r]_{s_1} & \cb{A_{0}^{1}}  \\
}
\end{align}
We know that:
\begin{itemize}
	\item $A_{0,0}^{1} \subset H_{1}$ is the braided Hopf algebra living in $\mathcal{YD}_{\C}(H_{0})$, which is obtained from the Hopf algebra projection $\big( d_{0} \colon H_{1} \to H_{0} \, , \, s_{0} \big)$.
	\item $A_{0,0}^{2} \subset H_{2}$ is the braided Hopf algebra living in $\mathcal{YD}_{\C}(H_{1})$, which is obtained from the Hopf algebra projection $\big( d_{0} \colon H_{2} \to H_{1} \, , \, s_{0} \big)$.
\end{itemize}

\begin{remark}
	We already know from \S \ref{interchange} that, $A_{0,0}^{1}$ also has a braided Hopf algebra structure in $\mathcal{YD}_{\C}(H_{1})$, considering the Hopf algebra projection $\big( d_{0} \colon H_{1} \to H_{0} \, , \, s_{0} \big)$.
\end{remark}

\begin{question}
	Considering the Hopf algebra morphisms \eqref{h2-case} given in $\C$, is it possible to obtain another induced Hopf algebra projection in $\yd$ between $A_{0,0}^{2}$ and $A_{0,0}^{1}$? 
\end{question}

\bigskip

\begin{Idea}
	We have two candidates for it, namely:
	$$\xymatrix@C=40pt@R=30pt{
		A_{0,0}^{2} \ar@<0.7ex>[r]^{d_1} & A_{0,0}^{1} \ar@<0.7ex>[l]^{s_1}} \, , \quad \quad \quad \xymatrix@C=40pt@R=30pt{
		A_{0,0}^{2} \ar@<0.7ex>[r]^{d_2} & A_{0,0}^{1} \ar@<0.7ex>[l]^{s_1} \, ,
		%	A_{0}^{2} \ar@<0.75ex>[r]^{d_1} \ar@<-0.75ex>@{<-}[r]_{s_1} & \cb{A_{0}^{1}}  \\
	}$$
	where the morphisms are the restrictions\footnote{We always use same notation for the restricted cases of face and degeneracy morphisms in the rest of the paper.}.
\end{Idea}

\begin{Lemma}\label{f-g-commutes}
	The following diagrams commute:
	$$\xymatrix@C=40pt@R=30pt{
		H_{2} \ar[r]^{d_2} \ar@<-0.5ex>[d]_{f^{2}_{0,0}} \ar@<0.5ex>[d]^{g^{2}_{0,0}} & H_{1} \ar@<0.5ex>[d]^{g^{1}_{0,0}} \ar@<-0.5ex>[d]_{f^{1}_{0,0}} \\
		A_{0,0}^{2} \ar[r]^{d_2} & A_{0,0}^{1} 
		%	A_{0}^{2} \ar@<0.75ex>[r]^{d_1} \ar@<-0.75ex>@{<-}[r]_{s_1} & \cb{A_{0}^{1}}  \\
	} \quad \quad \quad
	\xymatrix@C=40pt@R=30pt{
		H_{2}  \ar@<-0.5ex>[d]_{f^{2}_{0,0}} \ar@<0.5ex>[d]^{g^{2}_{0,0}} & H_{1} \ar@<0.5ex>[d]^{g^{1}_{0,0}} \ar@<-0.5ex>[d]_{f^{1}_{0,0}} \ar[l]_{s_1} \\
		A_{0,0}^{2} & A_{0,0}^{1} \ar[l]_{s_1}
		%	A_{0}^{2} \ar@<0.75ex>[r]^{d_1} \ar@<-0.75ex>@{<-}[r]_{s_1} & \cb{A_{0}^{1}}  \\
	} 
	$$
	In other words, the kernel generator maps $f^{\star}_{0,0}$ and $g^{\star}_{0,0}$ are compatible with $d_2$ and $i_1$ in $\C$. 
\end{Lemma}

\begin{Proof}
	For the $f^{\star}_{0,0}$ case, we have (by using simplicial identites):
	\begin{align*}
	f^{1}_{0,0} \, d_2 & = \nabla \, (\id \tn s_0 d_0 S) \, \Delta \, d_2 & f^{2}_{0,0} \, i_1& = \nabla \, (\id \tn s_0 d_0 S) \, \Delta \, s_1 \\
	& = \nabla \, (d_2 \tn s_0 d_0 S d_2) \, \Delta &&  = \nabla \, (s_1 \tn s_0 d_0 S s_1) \, \Delta \\
	& = \nabla \, (d_2 \tn s_0 d_0 d_2 S) \, \Delta && = \nabla \, (s_1 \tn s_0 d_0 s_1 S) \, \Delta \\
	& = \nabla \, (d_2 \tn s_0 d_1 d_0 S) \, \Delta && = \nabla \, (s_1 \tn s_0 s_0 d_0 S) \, \Delta \\
	& = \nabla \, (d_2 \tn d_2 s_0 d_0 S) \, \Delta && = \nabla \, (s_1 \tn s_1 s_0 d_0 S) \, \Delta \\
	& = d_2 \, \nabla (\id \tn s_0 d_0 S) \, \Delta && = s_1 \, \nabla (\id \tn s_0 d_0 S) \, \Delta \\
	& = d_2 \, f^{2}_{0,0} \, , && = s_1 \, f^{2}_{0,0} \, .
	\end{align*}
	%	\begin{align*}
	%		f^{1}_{0,0} \, d_2 & = \nabla \, (\id \tn i_0 d_0 S) \, \Delta \, d_2 \\
	%		& = \nabla \, (d_2 \tn i_0 d_0 S d_2) \, \Delta \\
	%		& = \nabla \, (d_2 \tn i_0 d_0 d_2 S) \, \Delta \\
	%		& = \nabla \, (d_2 \tn i_0 d_1 d_0 S) \, \Delta \\
	%		& = \nabla \, (d_2 \tn d_2 i_0 d_0 S) \, \Delta \\
	%		& = d_2 \, \nabla (\id \tn i_0 d_0 S) \, \Delta \\
	%		& = d_2 \, f^{2}_{0,0} \, ,
	%	\end{align*}
	%	and also:
	%	\begin{align*}
	%	f^{2}_{0,0} \, i_1& = \nabla \, (\id \tn i_0 d_0 S) \, \Delta \, i_1 \\
	%	& = \nabla \, (i_1 \tn i_0 d_0 S i_1) \, \Delta \\
	%	& = \nabla \, (i_1 \tn i_0 d_0 i_1 S) \, \Delta \\
	%	& = \nabla \, (i_1 \tn i_0 i_0 d_0 S) \, \Delta \\
	%	& = \nabla \, (i_1 \tn i_1 i_0 d_0 S) \, \Delta \\
	%	& = i_1 \, \nabla (\id \tn i_0 d_0 S) \, \Delta \\
	%	& = i_1 \, f^{2}_{0,0} \, .
	%	\end{align*}
	
	Same simplicial identities also proves the $g^{\star}_{0,0}$ case.
	%	
	%	\medskip
	%	
	%	 \centerline{\relabelbox
	%		\epsfysize 3cm
	%		\epsfbox{figs/f-g-commutes.eps}
	%		\relabel{1}{\scalebox{0.7}{$d_2$}}
	%		\relabel{3}{\scalebox{0.7}{$d_2$}}
	%		\relabel{2}{\scalebox{0.7}{$i_0$}}
	%		\relabel{4}{\scalebox{0.7}{$d_0$}}
	%		\relabel{5}{\scalebox{0.7}{$d_2$}}
	%		\relabel{8}{\scalebox{0.7}{$d_2$}}
	%		\relabel{6}{\scalebox{0.7}{$i_0$}}
	%		\relabel{7}{\scalebox{0.7}{$d_0$}}
	%		\relabel{9}{\scalebox{0.7}{$d_2$}}
	%		\relabel{10}{\scalebox{0.7}{$i_0$}}
	%		\relabel{11}{\scalebox{0.7}{$d_0$}}
	%		\relabel{12}{\scalebox{0.7}{$d_2$}}
	%		\relabel{13}{\scalebox{0.7}{$d_2$}}
	%		\relabel{14}{\scalebox{0.7}{$i_0$}}
	%		\relabel{15}{\scalebox{0.7}{$d_1$}}
	%		\relabel{16}{\scalebox{0.7}{$d_0$}}
	%		\relabel{20}{\scalebox{0.7}{$d_2$}}
	%		\relabel{17}{\scalebox{0.7}{$d_2$}}
	%		\relabel{18}{\scalebox{0.7}{$i_0$}}
	%		\relabel{19}{\scalebox{0.7}{$d_0$}}
	%		\relabel{21}{\scalebox{0.7}{$d_2$}}
	%		\relabel{22}{\scalebox{0.7}{$i_0$}}
	%		\relabel{23}{\scalebox{0.7}{$d_0$}}
	%		\relabel{24}{\scalebox{0.7}{$d_2$}}
	%%		\relabel{2}{$f^{\star} $}
	%		\relabel{e1}{$\overset{\scalebox{0.7}{$\mathrm{defn}$}}{=}$}
	%		\relabel{e7}{$\overset{\scalebox{0.7}{$\mathrm{defn}$}}{=}$}
	%		\relabel{e2}{$=$}
	%		\relabel{e6}{$=$}
	%%		\relabel{4}{$=$}
	%		\relabel{e3}{$\overset{\scalebox{0.7}{$\mathrm{(H5)}$}}{=}$}
	%		\relabel{e4}{$\overset{\scalebox{0.7}{$\mathrm{(1)}$}}{=}$}
	%			\relabel{e5}{$\overset{\scalebox{0.7}{$\mathrm{(3c)}$}}{=}$}
	%%		\relabel{6}{$=$}
	%%		\relabel{7}{$\overset{\scalebox{0.7}{$\mathrm{defn}$}}{=}$}
	%		\endrelabelbox}
	% \n Similarly, $g^{\star}$ commutes with both maps as well.
\end{Proof}

\begin{remark}
	However, $d_1$ is not compatible with $f^{\star}_{0,0}$ and $g^{\star}_{0,0}$.
\end{remark}

\begin{Lemma}\label{d-i-morphism}
	We have the following Hopf algebra projection in $\mathcal{YD}_{\C}(H_{1})$: 
	%	\begin{align}\label{new-projection}
	%	(d_{1} \colon A_{0}^{2} \to A_{0}^{1}, s_{1}) \, ,
	%	\end{align}
	\begin{align}\label{h2-case2}
	\xymatrix@C=40pt@R=30pt{
		A_{0,0}^{2} \ar@<0.7ex>[r]^{d_2} & A_{0,0}^{1} \ar@<0.7ex>[l]^{s_1} \, . \\
		%	A_{0}^{2} \ar@<0.75ex>[r]^{d_1} \ar@<-0.75ex>@{<-}[r]_{s_1} & \cb{A_{0}^{1}}  \\
	}
	\end{align}
\end{Lemma}

\begin{Proof}
	It is clear that $d_2,i_1$ are well-defined since the simplicial idendities $(1)$ and $(3a)$. They define algebra morphisms since $A_{0}^{2}$ and $A_{0}^{1}$ are sub algebras of $H_2$ and $H_1$, respectively. Moreover:
	\begin{itemize}
		\item $d_2,i_1$ are coalgebra (therefore bialgebra) morphisms, since:
		\begin{align*}
		\D^{}_{A_{0,0}^{1}} \, d_2 & = (f^{1}_{0,0} \tn \id) \, \D_{\C} \, d_2 \\
		& = (f^{1}_{0,0} \tn \id) \, (d_2 \tn d_2) \, \D_{\C} \\
		& = (d_2 \tn d_2) \, (f^{2}_{0,0} \tn \id) \, \D_{\C} \quad \mathrm{( Lemma \,\, \ref{f-g-commutes})} \\
		& = (d_2 \tn d_2) \, \D^{}_{A_{0,0}^{2}}
		\end{align*}
		that proves $d_2$ is a coalgebra morphism. Similarly, $s_1$ defines a coalgebra morphism as well.
		\item It is also straightforward from Lemma \ref{f-g-commutes} that, $d_2$ and $s_1$ are compatible with the antipode, by using the fact that $S^{\star} = g_{0,0}^{\star}$.
	\end{itemize}
	That proves $d_2$ and $s_1$ Hopf algebra morphisms in $\mathcal{YD}_{\C}(H_{1})$. Therefore, \eqref{h2-case2} defines a Hopf algebra projection in $\mathcal{YD}_{\C}(H_{1})$.
\end{Proof}
%
%\bigskip
% 
%\noindent We again have a Hopf algebra projection in category $\mathcal{DY}_{\C}(H_{1})$: $$(d_{1} \colon A_{0}^{2} \to A_{0}^{1}, i_{i})$$
%(Remark that: the maps are now restricted.)

\begin{remark}\label{problem1}
	Since $d_1$ does not commute with $f^{\star}_{0,0}$ and $g^{\star}_{0,0}$, it does not define a Hopf algebra morphism between $A_{0,0}^{2} \to A_{0,0}^{1}$ in $\mathcal{YD}_{\C}(H_{1})$. Therefore we can not obtain any Hopf algebra projetion through $d_1$.
\end{remark}

\begin{Corollary}
	Consequently, if we apply Radford's Theorem to \eqref{h2-case2} in  $\mathcal{YD}_{\C}(H_{1})$, we get the braided Hopf algebra: 
	\begin{align*}
	A_{2,1}^{2} = RKer_{\big( \mathcal{YD}_{\C}(H_{1}) \big)} \big( d_{2} \colon A_{0,0}^{2} \to A_{0,0}^{1} \big) \, ,
	\end{align*}
	which is living in the (induced) braided monoidal category $\mathcal{YD}_{\big( \mathcal{YD}_{\C}(H_{1}) \big)} (A_{0,0}^{1})$.
\end{Corollary}

\begin{remark}
	Note that, $A_{2,1}^{2} \subset A_{0,0}^{2}$ at the subalgebra level. 
	%in $\mathcal{YD}_{\C}(H_{1})$.
\end{remark}

%
%\bigskip
%
%which is subset (and sub algebra) of $A_{0}^{2}$.
%
%$$A_{1}^{2} \subset A_{0}^{2}$$ in the category $\mathcal{DY}_{\mathcal{DY}_{\C}(H_{1})}(A_{0}^{1})$.

\medskip

\begin{tcolorbox}
	\begin{remark}\label{summary1}
		We can summarize the above constructions as follows:
		\begin{align}\label{summary}
		\xymatrix@C=100pt@R=40pt{
			{H}_{2}\ar@<3ex>[r]|{d_2} \ar@<1.5ex>[r]|{d_1} \ar@<0ex>[r]|{d_0} &
			{H}_{1}\ar@<1.5ex>[r]|{d_1} \ar[r]|{d_0}
			\ar@<1.5ex>[l]|{s_0} \ar@<3ex>[l]|{s_1} &
			{H}_{0} \ar@<1.5ex>[l]|{s_0} \\
			A_{0,0}^{2} \ar@<0.75ex>[r]^{d_2} \ar@<-0.75ex>@{<-}[r]_{s_1} & A_{0,0}^{1} %\ar@<-0.5ex>@{}[ur]^(.1){}="a"^(.85){}="b" \ar@[red] "a";"b"_{\d_1}
			&  \\
			A_{2,1}^{2} %\ar@<-0.5ex>@{}[ur]^(.1){}="a"^(.85){}="b" \ar@[red] "a";"b"_{\d_1} & &
		}
		\end{align}
		with expressing that, each level consists of Hopf algebras (projections) in different categories.
	\end{remark}
\end{tcolorbox}

\medskip

\begin{Lemma}\label{twisted}
	Let us define $\partial _{1}\colon A_{0,0}^{1} \to H_0$ as the restriction of $d_1$. Then, $\d_1$ defines a twisted Hopf algebra map, i.e. the following diagram commutes:
	\begin{align*}
	\xymatrix@R=40pt@C=50pt{ A_{0,0}^{1}  \ar[r]^-{\d_1} \ar[d]_{\D_{A_{0,0}^{1}}} & H_0 \ar[r]^-{\D_{H_{0}}} & H_0 \tn H_0  \\ A_{0,0}^{1} \tn A_{0,0}^{1} \ar[rr]^-{\d_1 \tn \, \rho} & & H_0 \tn H_0 \tn A_{0,0}^{1} \ar[u]_-{\nabla \tn \d_1} } 
	\end{align*}
\end{Lemma}

\begin{Proof}
	By using the simplicial identities, we have:
	\begin{align*}
	( \nabla \tn d_1 ) (d_1 \tn \rho) \, \D_{A_{0,0}^{1}} (a) & = ( \nabla \tn d_1 ) (d_1 \tn \rho) (f_{0,0}^{1} \tn \id) \, \D_{H_0} (a) \\
	& = ( \nabla \tn d_1 ) (d_1 \tn \rho) (f_{0,0}^{1} \tn \id) \,\, \underset{(a)}{\sum } \, \, a' \tn a'' \\ 
	& = ( \nabla \tn d_1 ) \,\, \underset{(a)}{\sum } \,\, d_1 f_{0,0}^{1} (a') \tn \rho (a'')  \\
	& = ( \nabla \tn d_1 ) \,\, \underset{(a)}{\sum } \,\, d_1 f_{0,0}^{1} (a') \tn d_0 (a'') \tn a'''  \\
	& = ( \nabla \tn d_1 ) \,\, \underset{(a)}{\sum } \, \, d_1 (a') \, d_1 s_0 d_0 S (a'') \tn d_0 (a''') \tn a''''  \\
	& = ( \nabla \tn d_1 ) \,\, \underset{(a)}{\sum } \, \, d_1 (a') \, d_0 S (a'') \tn d_0 (a''') \tn a''''  \\
	& = ( \id \tn d_1 ) \,\, \underset{(a)}{\sum } \, \, d_1 (a') \, d_0 S (a'') \, d_0 (a''') \tn a''''  \\
	& = ( \id \tn d_1 ) \,\, \underset{(a)}{\sum } \, \,  d_1 (a')  \tn a''  \\
	& = \,\, \underset{(a)}{\sum } \, \, d_1 (a')  \tn d_1(a'') = \underset{(d_1(a))}{\sum } d_1 (a)'  \tn d_1(a)''  \\
	& = \D_{H_0} \, d_1 (a) \, ,
	\end{align*}
	for all $a \in A_{0,0}^{1}$, which proves that the diagram above commutes.
\end{Proof}

\begin{remark}
	If we similarly define $\partial _{1}\colon A_{2,1}^{2} \to A_{0,0}^{1}$ as the restriction of $d_1 \colon H_2 \to H_1$, then $\d_1$ defines only an algebra morphism, not a (twisted) Hopf algebra map.
\end{remark}

\subsection{A problem for higher dimensions}

Now, let us consider the first three parts of the simplicial Hopf algebra. Then, apply Radford's Theorem step by step analogously to \S \ref{dimension2} and see what happens.

\medskip

\begin{remark}
	In the diagram above, the restriction $s_2 \colon A_{2,1}^{2} \to A_{2,1}^{3}$ is not well-defined, therefore $s_2$ does not define a (Hopf) algebra morphism.
\end{remark}

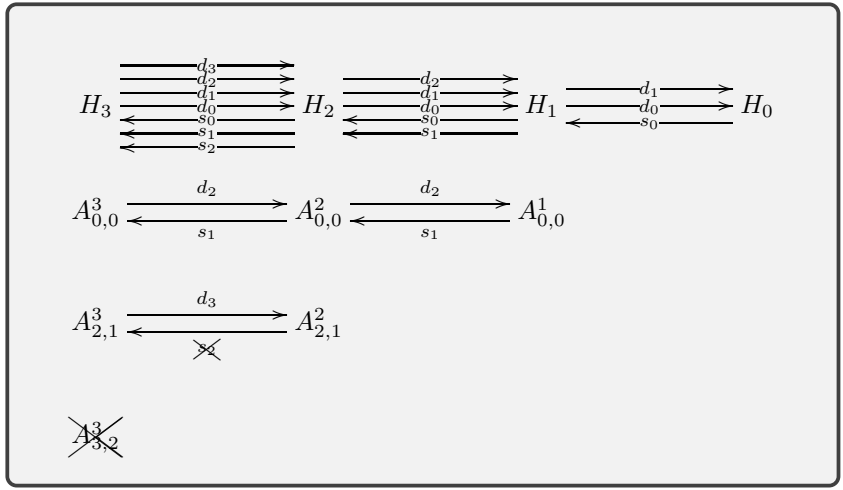
\begin{center}
\begin{tcolorbox}[text width = 10cm]
\begin{align*}
		\xymatrix@C=60pt@R=24pt{
		H_3  \ar@<3.6ex>[r]|{d_3} \ar@<2.4ex>[r]|{d_2} \ar@<1.2ex>[r]|{d_1} \ar@<0ex>[r]|{d_0}  &
		{H}_{2}\ar@<2.4ex>[r]|{d_2} \ar@<1.2ex>[r]|{d_1} \ar@<0ex>[r]|{d_0} \ar@<1.2ex>[l]|{s_0} \ar@<2.4ex>[l]|{s_1} \ar@<3.6ex>[l]|{s_2} &
		{H}_{1}\ar@<1.5ex>[r]|{d_1} \ar[r]|{d_0}
		\ar@<1.2ex>[l]|{s_0} \ar@<2.4ex>[l]|{s_1} &
		{H}_{0} \ar@<1.5ex>[l]|{s_0}  \\
		{ A_{0,0}^{3} } { \ar@<0.75ex>[r]^{d_2}}  
		{ \ar@<-0.75ex>@{<-}[r]_{s_1} }
		& { A_{0,0}^{2} } 
		 { \ar@<0.75ex>[r]^{d_2}}  
		{ \ar@<-0.75ex>@{<-}[r]_{s_1} } & 
		{ A_{0,0}^{1} } %\ar@<-0.5ex>@{}[ur]^(.1){}="a"^(.85){}="b" \ar@[red] "a";"b"_{\d_1}
		&  \\
		A_{2,1}^{3} { \ar@<0.75ex>[r]^{d_3} \ar@<-0.75ex>@{<-}[r]_{\xcancel{s_2} }} & A_{2,1}^{2}
		\\
		{
			\xcancel{A_{3,2}^{3} } }
		 & &
	}
\end{align*}
\end{tcolorbox}
\end{center}

Consequently, it is not possible to apply Radford's Theorem after the third level of the process since we do not have any candidate to obtain a Hopf algebra projection between $A_{2,1}^{3}$ and $A_{2,1}^{2}$ except $d_3 \, s_2$, in $\mathcal{YD}_{\big( \mathcal{YD}_{\C}(H_{2}) \big)} (A_{0,0}^{2})$.

\section{Application: Braided Hopf Crossed Modules}

In this section, our main is to understand braided Hopf crossed modules from the point of simplicial structures. But first, let us briefly recall the definition of braided Hopf crossed modules \cite{SM2}.

\subsection{Braided Hopf crossed modules}

\begin{Definition}
	Let $H$ be a Hopf algebra and $I$ be a braided Hopf algebra living in $\mathcal{YD}(H)$. The ``twisted Hopf algebra map'' $\d \colon I \to H$, namely an algebra morphism obeys: 
	$$ \D \big( \d(x) \big) = \sum_{(h)} \sum_{[h]} \d(x') \, x''_{H} \tn \d(x''_{I}) \, , $$ 
	is called a ``braided Hopf crossed module'' if it further satisfies:
	
	\begin{itemize}
		\item $\d(h \t x) = \sum_{(h)} h' \, \d(x) \, S(h'') ,$
		\item $ \d(x) \t y = \sum_{(x)} \sum_{[x]} x' \, \big( x''_{H} \t y \big) \,  S(x''_{I}) ,$
	\end{itemize}
	for all $x,y \in I$ and $h \in H$. Notice that, the right hand side of the second condition is the braided adjoint action in $\mathcal{YD}(H)$, given in \eqref{bad}. On the other hand, without the last condition, we call it a braided Hopf pre-crossed module.
\end{Definition}

\subsection{Generating the elements of $A^{2}_{2,1}$}

In this subsection, we will calculate a specific type of elements in a simplicial Hopf algebra $\mathcal{H}$ which is defined over the category of vector spaces.\footnote{Working in the category of vector spaces will allow us to have combinatorial calculations. On the other hand, our main aim is to make contact with Majid's braided crossed module notion \cite{SM2} where the base category is vector spaces. For the general case of such crossed modules, see \cite{zbMATH01036413}.} Afterwards, this specific type of elements will lead us to discover the relationship between simplicial Hopf algebra and Majid's braided Hopf crossed module definition. The following idea and terminology firstly defined and also applied in \cite{CC1,MP1} for the case of groups. They call it ``iterated Peiffer pairings" to construct higher level Peiffer elements to model higher dimensional categorical structures. However, in this paper, we slightly modify the idea to fit our construction.

\medskip

Let $\mathcal{H}$ be a simplicial Hopf algebra and recall the constructions given in \eqref{summary}. Then, for all $x,y \in A^{1}_{0,0}$, we have the following diagram to construct an element of $A^{2}_{2,1}$. 
%Remark that, the braided adjoint action in vector spaces is clearly normal adjoint action as seen at (\ref{adjoint}). Now let us to calculate Peiffer pairings of dimension $2$. 
\begin{align}\label{peiffer}
\xymatrix@R=20pt@C=20pt{
	A_{0,0}^{1} \times A_{0,0}^{1} \ar@{-->}[rr]^-{F_{(0)(1) }} \ar[dd]_{s_{0}\times s_{1}}
	&& A_{2,1}^{2} \\
	&& A_{0,0}^{2} \ar[u]_{f_{2,1}^{2}} \\
	H_{2} \times H_{2} \ar[rr]^{\vartriangleright_{ad}}
	&& H_{2} \ar[u]_{f_{0,0}^{2}}
} 
\end{align}

%Recall that from (\ref{rker-gen}) that: $$ A_{0}^{1} = RKer_{vect} \big( d_{0} \colon H_{1} \to H_{0} \big) $$
%For $n=2$, suppose $\alpha =(0)$, $\beta =(1)$ and $x,y\in
%NH_{1}=A_{0}^{1} $. Then:

\bigskip

If we calculate such elements, we get:
\begin{align*}
F_{(0)(1)}(x,y)& =f_{2,1}^{2}  \, f_{0,0}^{2} \, \Big( (s_{0}(x)\vartriangleright _{ad}s_{1}(y)
\Big) \\
& =f_{2,1}^{2}\left[ \underset{(s_{0}(x)\vartriangleright _{ad}s_{1}(x))}{\sum }%
\left[ s_{0}(x)\vartriangleright _{ad}s_{1}(y)\right] ^{\prime
} \,\, s_{0}d_{0}(S[s_{0}(x)\vartriangleright _{ad}s_{1}(y)]^{\prime \prime })\right] \\
& =f_{2,1}^{2}\left[ \underset{(x)(y)}{\sum } s_{0}(x') \, s_{1}(y') \, s_{0}S(x'''') \,\, s_{0}d_{0}S\Big(s_{0}(x'') \, s_{1}(y'') \, s_{0}S(x''') \Big) \right] \\
& =f_{2,1}^{2}\left[ \underset{(x)(y)}{\sum } s_{0}(x') \, s_{1}(y') \, s_{0}S(x'''') \,\, S\Big(s_{0}(x'') \, s_{0}d_{0}s_{1}(y'') \, s_{0}S(x''') \Big) \right] \\
& =f_{2,1}^{2}\left[ \underset{(x)(y)}{\sum } s_{0}(x') \, s_{1}(y') \, s_{0}S(x'''') \,\, S\Big(s_{0}(x'') \, s_{0}s_{0}d_{0}(y'') \, s_{0}S(x''') \Big) \right] \\
& =f_{2,1}^{2}\left[ \underset{(x)}{\sum } s_{0}(x') \, s_{1}(y) \, s_{0}S(x'''') \,\, S\Big(s_{0}(x'') \, s_{0}S(x''') \Big) \right] \quad \textrm{since} \,\, \Delta(y) \subset H_{1} \tn RKerd(d_{0}) \\
& =f_{2,1}^{2}\left[ \underset{(x)}{\sum } s_{0}(x') \, s_{1}(y) \, s_{0}S(x'''') \,\, S\Big(s_{0}(x'' \, S(x''')) \Big) \right] \\
& =f_{2,1}^{2}\left[ \underset{(x)}{\sum } s_{0}(x') \, s_{1}(y) \, s_{0}S(x''') \,\, S\Big( \epsilon(x'') \Big) \right] \\
& =f_{2,1}^{2}\left[ \underset{(x)}{\sum } s_{0}(x') \, s_{1}(y) \, s_{0}S(x'') \right] \\
& =f_{2,1}^{2} \Big( s_{0}(x)\vartriangleright _{ad}s_{1}(y) \Big)
\end{align*}
Here, we see that the element $\Big( s_{0}(x)\vartriangleright _{ad}s_{1}(y) \Big)$
is also in $A_{0,0}^{2}$. Now we need to calculate the element:
\begin{align}\label{element1}
f_{2,1}^{2} \Big( s_{0}(x)\vartriangleright _{ad}s_{1}(y) \Big).
\end{align}

\begin{remark}\label{remark1}
To continue the calculation, we need to see what $f_{2,1}^{2} $ is. We know from \eqref{aa} that:
\[
f_{2,1}^{2} = \nabla \, \big( id \tn s_{1}d_{2}S_{A_{0,0}^{2}} \big) \, \D_{A_{0,0}^{2}} \, ,
\]
where:
\begin{align*}
\D_{A_{0,0}^{2}} = \big( f_{0,0}^{2} \tn id \big) \, \D \, ,
\end{align*}
from \eqref{bb}. Therefore, in an explicit formula:
\begin{align*}
f_{2,1}^{2} (a) & = \nabla \, \big( id \tn s_{1}d_{2}S_{A_{0,0}^{2}} \big) \big( f_{0,0}^{2} \tn id \big) \, \D  \\
& = \nabla \, \big( id \tn s_{1}d_{2}S_{A_{0,0}^{2}} \big) \, \underset{(a)}{\sum } f_{0,0}^{2}(a') \tn a'' \\
& = \underset{(a)}{\sum } \, f_{0,0}^{2}(a') \,\, s_{1}d_{2}g_{0,0}^{2}(a'') \\
& = \underset{(a)}{\sum } \, a' \,\, s_{0}d_{0}S(a'') \,\, s_{1}d_{2} \Big( s_{0}d_{0}(a''') \, S(a'''') \Big) \\
& = \underset{(a)}{\sum } \, a' \,\, s_{0}d_{0}S(a'') \,\, s_{1}d_{2}s_{0}d_{0}(a''') \,\, s_{1}d_{2}S(a'''') \, .
\end{align*}
\end{remark}

\begin{remark}\label{remark2}
We also need to obtain the clear formulae of: $$\Big( \Delta \circ \Delta \Big) \Big( s_{0}(x) \ad s_{1}(y) \Big).$$

Since:
\begin{align*}
\Delta(x \ad y) = \underset{(x)(y)}{\sum } x' \, y' \, S(x'''') \tn x'' \, y'' \, S(x''') \, ,
\end{align*}
we get:
\begin{align*}
\Delta \Big( \Delta(x \ad y) \Big) = \underset{(x)(y)}{\sum } x' \, y' \, S(x^{viii}) \tn x'' \, y'' \, S(x^{vii}) \tn x''' \, y''' \, S(x^{vi}) \tn x'''' \, y'''' \, S(x^{v})
\end{align*}
Therefore:
\begin{align*}
& \Big( \Delta \circ \Delta \Big) \Big( s_{0}(x) \ad s_{1}(y) \Big) \\
& = \underset{(x)(y)}{\sum } s_{0}(x') \, s_{1}(y') \, s_{0}S(x^{viii}) \tn s_{0}(x'') \, s_{1}(y'') \, s_{0}S(x^{vii}) \tn s_{0}(x''') \, s_{1}(y''') \, s_{0}S(x^{vi}) \tn s_{0}(x'''') \, s_{1}(y'''') \, s_{0}S(x^{v}) \, .
\end{align*}
\end{remark}

\medskip

If we put $a = s_{0}(x) \ad s_{1}(y)$ in Remark \ref{remark1} via Remark \ref{remark2}, we therefore obtain the result of \eqref{element1} as follows: 
\begin{align*}
f_{1}^{2}(a)& =f_{1}^{2} \big(s_{0}(x)  \ad s_{1}(y) \big) \\ & = \underset{(x)(y)}{\sum } s_{0}(x') \, s_{1}(y') \, s_{0}S(x^{viii}) \,\, s_{0}d_{0}S \Big( s_{0}(x'') \, s_{1}(y'') \, s_{0}S(x^{vii}) \Big) \\
& \q\q\q\q\q\q \cb{s_{1}d_{2}s_{0}d_{0}} \Big( s_{0}(x''') \, s_{1}(y''') \, s_{0}S(x^{vi}) \Big) \,\, \cb{s_{1}d_{2}}S \Big( s_{0}(x'''') \, s_{1}(y'''') \, s_{0}S(x^{v}) \Big) \\
\\
& = \underset{(x)(y)}{\sum } s_{0}(x') \, s_{1}(y') \, s_{0}S(x^{viii}) \,\, s_{0}S(S(x^{vii})) \, s_{0}d_{0}s_{1}S(y'') \, s_{0}S(x'') \\ 
& \q\q\q\q\q\q \cb{s_{1}d_{2}s_{0}}(x''') \, \cb{s_{1}d_{2}s_{0}d_{0}s_{1}}(y''') \, \cb{s_{1}d_{2}s_{0}}S(x^{vi}) \,\, \cb{s_{1}d_{2}s_{0}}S(S(x^{v})) \, s_{1}S(y'''') \, \cb{s_{1}d_{2}s_{0}} S(x'''') \\
%\\
%& = \underset{(x)(y)}{\sum } s_{0}(x') \, s_{1}(y') \, \cb{s_{0}S(x^{viii}) \,\, s_{0}S(S(x^{vii}))} \, s_{0}d_{0}s_{1}S(y'') \, s_{0}S(x'') \\
%& \q\q\q\q\q\q s_{1}(x''') \, s_{1}d_{0}s_{1}(y''') \, \cb{s_{1}S(x^{vi}) \,\, s_{1}S(S(x^{v}))} \, s_{1}S(y'''') \, s_{1}S(x'''') \\
\\
& = \underset{(x)(y)}{\sum } s_{0}(x') \, s_{1}(y') \, s_{0}d_{0}s_{1}S(y'') \, s_{0}S(x'') \, \cb{s_{1}d_{2}s_{0}}(x''') \, \cb{s_{1}d_{2}s_{0}d_{0}s_{1}} (y''') \, s_{1}S(y'''') \,  \cb{s_{1}d_{2}s_{0}} S(x'''') \\
\\
& = \underset{(x)(y)}{\sum } s_{0}(x') \, s_{1}(y') \, s_{0}d_{0}s_{1}S(y'') \, s_{0}S(x'') \, \cb{s_{1}d_{2}s_{0}}(x''') \, \cb{s_{1}d_{0}s_{1}} (y''') \, s_{1}S(y'''') \,  \cb{s_{1}d_{2}s_{0}} S(x'''') \\
\end{align*}
by using the fact that $s_{1}d_{2}s_{0}d_{0}s_{1} = s_{1}d_{2}d_{0}s_{1}s_{1} = s_{1}d_{0}d_{3}s_{1}s_{1} = s_{1}d_{0}d_{3}s_{2}s_{1} = s_{1}d_{0}s_{1}$. Moreover, we also used:
\begin{align*}
\nabla \Big( \underset{(S(x))}{\sum} \, \big( S(x) \big)' \tn S\big( S(x) \big)'' \Big) = \nabla \Big( \underset{(x)}{\sum} \, S(x'') \tn S(S(x')) \Big) = \epsilon (x) \, 1 \, .
\end{align*}
%which obtained from:
%\[
%\nabla \Big( \underset{(x)}{\sum} \, x' \tn S(x'') \Big) = \epsilon (x)
%\]

\begin{Corollary}\label{peifferelement}
	If $x,y \in A_{0,0}^{1}$, then the element:
	\begin{align*}
	F_{(0)(1)}(x,y) = \underset{(x)(y)}{\sum } s_{0}(x') \, s_{1}(y') \, s_{0}d_{0}s_{1}S(y'') \, s_{0}S(x'') \, \cb{s_{1}d_{2}s_{0}}(x''') \, \cb{s_{1}d_{0}s_{1}} (y''') \, s_{1}S(y'''') \,  \cb{s_{1}d_{2}s_{0}} S(x'''')
	\end{align*}
	belongs to $A_{2,1}^{2}$, obtained from the diagram \eqref{peiffer}.
\end{Corollary}

\subsection{Braided Hopf crossed module through a simplicial structure}

We know from Lemma \ref{twisted} that, for a given simplicial Hopf algebra $\mathcal{H}$, there exists a twisted Hopf algebra map:
\begin{align}\label{xmod}
\d_1 \colon A_{0,0}^{1} \longrightarrow H_0 \, ,
\end{align}
in the sense of \eqref{summary1}.

\begin{Theorem}\label{maintheorem}
	Let $\mathcal{H'}$ be a simplicial Hopf algebra such that $A_{2,1}^{2}$ is the zero object. Then \eqref{xmod} gives rise to a braided Hopf crossed module structure where the action of $H_0$ on $A_{0,0}^{1}$ is defined by the adjoint action via degeneracy morphism $s_1$, namely:
	\begin{align*}
	h \vartriangleright x = s_{0} (h) \ad x = \underset{(h)}{\sum } s_{0}(h') \, x \, s_{0} \big( S(h'') \big) \, ,
	\end{align*}
	for all $h \in H_0$ and $x \in A_{0,0}^{1}$.
	
\end{Theorem}

\begin{Proof}
	
	\medskip
	
	$\spadesuit$ Clearly, we have a ``braided Hopf pre-crossed module'' with:
	\begin{align*}
	\partial_1 (h \vartriangleright x) & = d_{1} \Big( \underset{(h)}{\sum } s_{0}(h') \, x \, s_{0} \big( S(h'') \big) \Big) \\
	& = \underset{(h)}{\sum } d_{1} s_{0}(h') \, d_{1} (x) \, d_{1}s_{0} \big( S(h'') \big) \\
	& = \underset{(h)}{\sum } h' \, \d (x) \, S(h'') \, ,
	\end{align*}
	by using simplicial identities.
	
	\bigskip
	
	$\spadesuit$ To make it a ``braided Hopf crossed module'', we need to show that:
	\begin{align*}
	\d_1 (x) \vartriangleright y = \underset{[x]}{\sum } \u{x} \, (\u {\u x}_{H} \vartriangleright y) \, \u{S}(\u {\u x}_{A}) \, .
	\end{align*}
	%where the right hand side is called the braided adjoint action in $\mathcal{YD}(H)$.
	
	\medskip
	
	\noindent The left hand side is:
	\begin{align*}
	\d_1 (x) \vartriangleright y & = d_{1}(x) \vartriangleright y \\
	& = s_{0}d_{1}(x) \ad y \\
	& = \underset{(x)}{\sum } s_{0}d_{1}(x') \,\, y \,\, s_{0} \big( S(d_{1}(x'')) \big) \, ,
	\end{align*}
	while on the right hand side, we have: 
	
	\begin{align*}
	\underset{[x]}{\sum } \u{x} \, (\u {\u x}_{H} \vartriangleright y) \, \u{S}(\u {\u x}_{A}) & = \underset{(x)}{\sum } f(x') \, (x''_{H} \vartriangleright y) \, \u{S}(x''_{A}) \\
	& = \underset{(x)}{\sum } f(x') \, (d_{0}(x'') \vartriangleright y) \, \u{S}(x''') \\
	& = \underset{(x)}{\sum } f(x') \, \Big(s_{0} d_{0}(x'') \ad y \Big) \, \u{S}(x''') \\
	& = \underset{(x)}{\sum } (x')' \, S s_{0} d_{0}(x')'' \, \Big(s_{0} d_{0}(x'') \ad y \Big) \, \u{S}(x'') \\
	& = \underset{(x)}{\sum } x' \, S s_{0} d_{0}(x'') \, \Big(s_{0} d_{0}(x''') \ad y \Big) \, \u{S}(x'''') \\
	& = \underset{(x)}{\sum } x' \, S s_{0} d_{0}(x'') \, s_{0} d_{0}(x''') \, y \, S s_{0} d_{0} (x'''') \, \u{S}(x''''') \\
	& = \underset{(x)}{\sum } x' \, \e(x'').1 \,y \, S s_{0} d_{0} (x''') \, \u{S}(x'''') \\
	& = \underset{(x)}{\sum } x' \, y \, S s_{0} d_{0} (x'') \, \u{S}(x''') \\
	& = \underset{(x)}{\sum } x' \, y \, S s_{0} d_{0} (x'') \, g(x''') \\
	& = \underset{(x)}{\sum } x' \, y \, S s_{0} d_{0} (x'') \, s_{0}d_{0}(x''') \, S(x'''') \\
	& = \underset{(x)}{\sum } x' \, y \, \epsilon(x'') \, S(x''') \\
	& = \underset{(x)}{\sum } x' \, y \, S(x'') \\
	& = x \ad y
	\end{align*}
	which means braided adjoint action is equal to normal adjoint action in this case\footnote{Of course, this is not true in general. In fact, this is the consequence of simplicial identities and the definition of action and coaction.}.
	
	\medskip
	
	So, we need to show that:
	\begin{align}\label{main}
	\underset{(x)}{\sum } s_{0}d_{1}(x') \,\, y \,\, s_{0} \big( S(d_{1}(x'')) \big) = x \ad y \, .
	\end{align}
	
	Recall from Corollary \ref{peifferelement} that,
	\[
	t = \underset{(x)(y)}{\sum } s_{0}(x') \, s_{1}(y') \, s_{0}d_{0}s_{1}S(y'') \, s_{0}S(x'') \, \cb{s_{1}d_{2}s_{0}}(x''') \, \cb{s_{1}d_{0}s_{1}} (y''') \, s_{1}S(y'''') \,  \cb{s_{1}d_{2}s_{0}} S(x'''') \, ,
	\]
	belongs to $A_{2,1}^{2}$. 
	
	\medskip
	
	On the other hand, we know that $A_{2,1}^{2}$ is fixed and trivial (i.e. zero object), consequently we have:
	\[
	d_{1}(t) = \epsilon(t) 1_{H_1} \, .
	\]
	Therefore:
	\begin{align*}
	d_{1}(t) & = d_{1} \left( \underset{(x)(y)}{\sum } s_{0}(x') \, s_{1}(y') \, s_{0}d_{0}s_{1}S(y'') \, s_{0}S(x'') \, \cb{s_{1}d_{2}s_{0}}(x''') \, \cb{s_{1}d_{0}s_{1}} (y''') \, s_{1}S(y'''') \,  \cb{s_{1}d_{2}s_{0}} S(x'''')  \right) \\
	& = \underset{(x)(y)}{\sum } x' \, y' \, d_{0}s_{1}S(y'') \, S(x'') \, \cb{d_{2}s_{0}}(x''') \, \cb{d_{0}s_{1}} (y''') \, S(y'''') \,  \cb{d_{2}s_{0}} S(x'''')
	\end{align*}
	is equal to:
	\begin{align*}
	\epsilon(x)\epsilon(y)1_{H_{1}}
	\end{align*}
	
	Let $f \colon H\otimes H \to H$ be:
	\begin{align*}
	f(x\otimes y)= \underset{(x)(y)}{\sum } x' \, y' \, d_{0}s_{1}S(y'') \, S(x'') \, \cb{d_{2}s_{0}}(x''') \, \cb{d_{0}s_{1}} (y''') \, S(y'''') \,  \cb{d_{2}s_{0}} S(x'''') \, ,
	\end{align*}
	and therefore, we have:
	\begin{align*}
	f(x\otimes y)=\epsilon(x)\epsilon(y)1_{H_{1}} \, .
	\end{align*}
	
	We can easily write,
	\begin{align*}
	 \underset{(x)(y)}{\sum } f(x'\otimes y')\otimes x''\otimes y''=  \underset{(x)(y)}{\sum } \left(\epsilon(x')\epsilon(y')1_{H}\right)\otimes x''\otimes y'' \, .
	\end{align*}
	
	Moreover, by using the diagram:
	$$\xymatrix@R=40pt@C=20pt{
		H \otimes H \ar[rrr]^{\Delta \otimes \Delta} \ar@{.>}[dd] &&&
		H \otimes H \otimes H \otimes H \ar[d]^{\id \otimes \tau \otimes \id} \\
		&&& H \otimes H \otimes H \otimes H \ar@/^1.2pc/[d]^{\nabla(\epsilon \otimes \epsilon) \otimes d_2 s_0 \otimes \id}
		\ar@/_1.2pc/[d]_{f \otimes d_2 s_0 \otimes \id} \\
		H \otimes H
		&&& H \ar[lll]^{m(\id \otimes m)}
	}$$
	we can obtain:	
	\begin{align}\label{p1}
	 \underset{(x)(y)}{\sum } f(x'\otimes y') \, \big( d_2 s_0 (x'') \, y'' \big) =  \underset{(x)(y)}{\sum } \left(\epsilon(x')\epsilon(y')1_{H}\right) \, \big( d_2 s_0 (x'') \, y'' \big)
	\end{align}
	The right hand side is equal to:
	\[
	d_2 s_0(x) \, y
	\]
	while on the left hand side (after the calculations) we have:
	\begin{align}\label{calc1}
	\begin{split}
 	\underset{(x)(y)}{\sum } f(x'\otimes y') \,& \big( d_2 s_0 (x'') \, y'' \big) \\ 
 	& = \underset{(x)(y)}{\sum } x' \, y' \, d_{0}s_{1}S(y'') \, S(x'') \, \cb{d_{2}s_{0}}(x''') \, \cb{d_{0}s_{1}} (y''') \, S(y'''') \,  \cb{d_{2}s_{0}} S(x'''') \, \cb{d_{2}s_{0}} (x''''') y'''''  \\
	& = \underset{(x)(y)}{\sum } x' \, y' \, d_{0}s_{1}S(y'') \, S(x'') \, \cb{d_{2}s_{0}}(x''') \, \cb{d_{0}s_{1}} (y''') \, .
	\end{split}
	\end{align}
	Since $y \in RKer(d_{0})$, we have $$\underset{(y)}{\sum } \, y' \tn y'' \tn y''' \subseteq H_1 \tn H_1 \tn  RKer(d_{0}) \, .$$ 
	Therefore continuing with \eqref{calc1}, we have:
	\begin{align*}
	& =  \underset{(x)(y)}{\sum } x' \, y' \, d_{0}s_{1}S(y'') \, S(x'') \, \cb{d_{2}s_{0}}(x''') \, \cb{d_{0}s_{1}} (y''') \\
	& =  \underset{(x)(y)}{\sum } x' \, y' \, s_{0}d_{0}S(y'') \, S(x'') \, \cb{d_{2}s_{0}}(x''') \, \cb{s_{0}d_{0}} (y''') \\
	& = \underset{(x)(y)}{\sum } x' \, y' \, s_{0}d_{0}S(y'') \, S(x'') \, \cb{d_{2}s_{0}}(x''') 
	\end{align*}
	Again, since $y \in RKer(d_{0})$, we have $\underset{(y)}{\sum } \, y' \tn y'' \subseteq H_1 \tn RKer(d_{0}) $. Consequently:
	\begin{align*}
	& = \underset{(x)(y)}{\sum } x' \, y' \, s_{0}Sd_{0}(y'') \, S(x'') \, \cb{d_{2}s_{0}}(x''') \\
	& = \underset{(x)}{\sum } x' \, y \, S(x'') \, \cb{d_{2}s_{0}}(x''') \, .
	\end{align*}
	
	So we have a new equality:
	\begin{align*}
	\underset{(x)}{\sum } x' \, y \, S(x'') \, \cb{d_{2}s_{0}}(x''') = d_2 s_0(x) \, y
	\end{align*}
	
	With the similar idea above, let $g \colon H\otimes H \to H$ be:
	\begin{align*}
	g(x\otimes y)= \underset{(x)}{\sum } x' \, y \, S(x'') \, \cb{d_{2}s_{0}}(x''')
	\end{align*}
	Therefore:
	\begin{align*}
	g(x\otimes y) = d_2 s_0(x) \, y
	\end{align*}
	that yields:
	\begin{align*}
	 \underset{(x)}{\sum } \, g(x'\otimes y)\otimes x'' =  \underset{(x)}{\sum } \, \big( d_2 s_0(x') \, y \big) \otimes x''
	\end{align*}
%	and also:
%	\begin{align*}
%	g(x'\otimes y')\otimes S(x'') \otimes y''= \left( x' \, y' \right) \otimes S(x'') \otimes y''
%	\end{align*}
%	And therefore:
%	\begin{align*}
%	g(x'\otimes y') \,\, S(x'') \,  \epsilon(y'') = \left( x' \, y' \right) \,\, S(x'') \, \epsilon (y'')
%	\end{align*}
	Again, consider the following diagram:
	$$\xymatrix@R=40pt@C=20pt{
		H \otimes H \ar[rrr]^{\Delta \otimes \id} \ar@{.>}[dd] &&&
		H \otimes H \otimes H  \ar[d]^{\id \otimes \tau} \\
		&&& H \otimes H \otimes H \ar@/^1.2pc/[d]^{\nabla(d_2s_0 \otimes \id) \otimes d_2 s_0 S}
		\ar@/_1.2pc/[d]_{g \otimes d_2 s_0 S} \\
		H
		&&& H \otimes H \ar[lll]^{\nabla}
	}$$

	By using simplicial identities, on the right hand we have:
	\begin{align*}
	\underset{(x)}{\sum } s_{0}d_{1}(x') \,\, y \,\, s_{0} \big( S(d_{1}(x'')) \big) \, ,
	\end{align*}
	and on the left hand side we get:
	\begin{align*}
	g(x'\otimes y') \,\, d_2 s_0 S(x'') & = x \ad y \, ,
	\end{align*}
	which satisfies the eqauality (\ref{main}) and gives us the second condition of braided Hopf crossed modules.
\end{Proof}

\section{Conclusion}

The Moore complex \cite{Moore1} of a simplicial group is a chain complex:
\begin{align}\label{moore1}
N(\mathcal{G})= \big( \dots \ra{d_{(n+1)}}  N(G)_n \ra{d_n} \dots  \ra{d_{3}} N(G)_2 \ra{d_2} N(G)_1 \ra{d_1} G_0  \big)
\end{align}
of groups, where $N(G)_n=\bigcap_{i=0}^{n-1} \ker(d_i)$ at level $n$, and the boundary morphisms $d_n\colon N(G)_n \to N(G)_{n-1}$ are the restrictions of $d_n\colon G_n \to G_{n-1}$. Moreover, $N(\mathcal{G})$ defines a normal chain complex of groups, namely $d_n(N(G)_n) \trianglelefteq N(G)_{n-1}$, for all $n \geq 1$. Thus, the Moore complex can be considered as the normalized chain complex of a simplicial group. 

\medskip

Based on this definition, we know that, the category of simplicial groups whose Moore complex is with length one, is equivalent to the category of group crossed modules \cite{Menagerie} (for the monoid version, see \cite{Bohm3}). The proof (for one direction) briefly contains the following functor: In such a simplicial group, the 2-truncation of \eqref{moore1}, namely:  
\begin{align}\label{moore0}
	\underset{0}{\underbrace{\ker(d_0) \cap \ker(d_1)}} \ra{d_2} \ker(d_0) \ra{d_1} G_0 \, ,
\end{align} 
defines a group crossed module. The proof mainly uses the specific type of elements of $\ker(d_0) \cap \ker(d_1)$ which is already trivial from the assumption. 

\medskip

However, if we slightly modify the definition of the Moore complex as:
\begin{itemize}
\item $N(G)_n=\bigcap_{i=0}^{n} \ker(d_i)$ where $i \neq 1$, 
\item put $\d_n \colon N(G)_n \to N(G)_{n-1}$ as the restrictions of  $d_1\colon G_n \to G_{n-1}$,
\end{itemize}
then, we obtain an alternative definition of \eqref{moore1} (through simplicial identities). Moreover, we already checked that this alternative definition still involves a group crossed module analogously to \eqref{moore0} as follows:
\begin{align}\label{moore2}
	 \underset{0}{\underbrace{\ker(d_0) \cap \ker(d_2)}} \ra{d_1} \ker(d_0) \ra{d_1} G_0 \, ,
\end{align}
after similar type of calculations. This gives an alternative functor from the category of that of simplicial groups to the category of group crossed modules.

\medskip

Now, let us go back to the category of Hopf algebras. In fact, Theorem \ref{maintheorem} generalizes the group theoretical case \eqref{moore2} to the category of Hopf algebras. However, we can not obtain an anologous version of \eqref{moore0} in the category of Hopf algebras due to Remark \ref{problem1} which was one of the major problems we faced during our research. From this point of view, this paper can also be considered as the first serious approach to understand the Moore complex of simplicial Hopf algebras in the most general case.

\bibliographystyle{abbrv}
\addcontentsline{toc}{section}{References}
\bibliography{Moore2}

\begin{thebibliography}{10}

\bibitem{Agore}
A.~L. {Agore}.
\newblock {Categorical constructions for Hopf algebras.}
\newblock {\em {Commun. Algebra}}, 39(4):1476--1481, 2011.

\bibitem{AgoreII}
A.~L. {Agore}.
\newblock {Limits of coalgebras, bialgebras and Hopf algebras.}
\newblock {\em {Proc. Am. Math. Soc.}}, 139(3):855--863, 2011.

\bibitem{AD1}
N.~{Andruskiewitsch} and J.~{Devoto}.
\newblock {Extensions of Hopf algebras.}
\newblock {\em {St. Petersbg. Math. J.}}, 7(1):22--61, 1995.

\bibitem{zbMATH01824088}
N.~{Andruskiewitsch} and H.-J. {Schneider}.
\newblock {Pointed Hopf algebras.}
\newblock In {\em {New directions in Hopf algebras}}, pages 1--68. Cambridge:
  Cambridge University Press, 2002.

\bibitem{zbMATH01036413}
Y.~N. {Bespalov}.
\newblock {Crossed modules and quantum groups in braided categories.}
\newblock {\em {Appl. Categ. Struct.}}, 5(2):155--204, 1997.

\bibitem{zbMATH07140994}
G.~{B\"ohm}.
\newblock {Crossed modules of monoids. I: Relative categories.}
\newblock {\em {Appl. Categ. Struct.}}, 27(6):641--662, 2019.

\bibitem{Bohm2}
G.~{B\"ohm}.
\newblock {Crossed modules of monoids. II: Relative crossed modules.}
\newblock {\em {Appl. Categ. Struct.}}, 2020.

\bibitem{Bohm3}
G.~{B\"ohm}.
\newblock {Crossed modules of monoids. III: Simplicial monoids of Moore length
1.}
\newblock arXiv:1803.04622.

\bibitem{zbMATH06456664}
Y.~{Boyaci}, J.~M. {Casas}, T.~{Datuashvili}, and E.~O. {Uslu}.
\newblock {Actions in modified categories of interest with application to
  crossed modules.}
\newblock {\em {Theory Appl. Categ.}}, 30:882--908, 2015.

\bibitem{zbMATH06825034}
R.~{Brown}.
\newblock {Modelling and computing homotopy types: I.}
\newblock {\em {Indag. Math., New Ser.}}, 29(1):459--482, 2018.

\bibitem{zbMATH01944940}
R.~{Brown} and I.~{\.I\c{c}en}.
\newblock {Homotopies and automorphisms of crossed modules of groupoids.}
\newblock {\em {Appl. Categ. Struct.}}, 11(2):185--206, 2003.

\bibitem{zbMATH03521203}
R.~{Brown} and C.~{Spencer}.
\newblock {$\mathcal G$-groupoids, crossed modules and the fundamental groupoid
  of a topological group.}
\newblock {\em {Nederl. Akad. Wet., Proc., Ser. A}}, 79:296--302, 1976.

\bibitem{zbMATH07003760}
D.~{Bulacu}, S.~{Caenepeel}, F.~{Panaite}, and F.~{Van Oystaeyen}.
\newblock {\em {Quasi-Hopf algebras. A categorical approach.}}, volume 171.
\newblock Cambridge: Cambridge University Press, 2019.

\bibitem{SB2}
S.~{Burciu}.
\newblock {Categorical Hopf kernels and representations of semisimple Hopf
  algebras.}
\newblock {\em {J. Algebra}}, 337(1):253--260, 2011.

\bibitem{zbMATH06005967}
S.~{Burciu}.
\newblock {Kernels of representations and coideal subalgebras of Hopf
  algebras.}
\newblock {\em {Glasg. Math. J.}}, 54(1):107--119, 2012.

\bibitem{CC1}
P.~{Carrasco} and A.~M. {Cegarra}.
\newblock {Group-theoretic algebraic models for homotopy types.}
\newblock {\em {J. Pure Appl. Algebra}}, 75(3):195--235, 1991.

\bibitem{E1}
G.~J. {Ellis}.
\newblock {Homotopical aspects of Lie algebras.}
\newblock {\em {J. Aust. Math. Soc., Ser. A}}, 54(3):393--419, 1993.

\bibitem{1905.09620}
K.~Emir.
\newblock {The Moore complex of a simplicial cocommutative Hopf algebra.}
\newblock arXiv:1905.09620.

\bibitem{JFM1}
J.~{Faria Martins}.
\newblock {Crossed modules of Hopf algebras and of associative algebras and
  two-dimensional holonomy.}
\newblock {\em {J. Geom. Phys.}}, 99:68--110, 2016.

\bibitem{zbMATH05815821}
J.~{Faria Martins} and R.~{Picken}.
\newblock {On two-dimensional holonomy.}
\newblock {\em {Trans. Am. Math. Soc.}}, 362(11):5657--5695, 2010.

\bibitem{zbMATH07075977}
M.~{Gran}, F.~{Sterck}, and J.~{Vercruysse}.
\newblock {A semi-abelian extension of a theorem by Takeuchi.}
\newblock {\em {J. Pure Appl. Algebra}}, 223(10):4171--4190, 2019.

\bibitem{zbMATH02096926}
G.~{Janelidze}.
\newblock {Internal crossed modules.}
\newblock {\em {Georgian Math. J.}}, 10(1):99--114, 2003.

\bibitem{JS2}
A.~Joyal and R.~Street.
\newblock Braided tensor categories.
\newblock {\em Advances in Mathematics}, 102(1):20--78, 11 1993.

\bibitem{zbMATH06812675}
V.~{Lebed} and F.~{Wagemann}.
\newblock {Representations of crossed modules and other generalized
  Yetter-Drinfel'd modules.}
\newblock {\em {Appl. Categ. Struct.}}, 25(4):455--488, 2017.

\bibitem{SM2}
S.~Majid.
\newblock Strict quantum 2-groups.
\newblock arXiv:1208.6265.

\bibitem{zbMATH00487164}
S.~{Majid}.
\newblock {Braided matrix structure of the Sklyanin algebra and of the quantum
  Lorentz group.}
\newblock {\em {Commun. Math. Phys.}}, 156(3):607--638, 1993.

\bibitem{SM3}
S.~{Majid}.
\newblock {Algebras and Hopf algebras in braided categories.}
\newblock In {\em {Advances in Hopf algebras. Conference, August 10-14, 1992,
  Chicago, IL, USA}}, pages 55--105. New York, NY: Marcel Dekker, 1994.

\bibitem{SM1}
S.~{Majid}.
\newblock {\em {Foundations of quantum group theory.}}
\newblock Cambridge: Cambridge Univ. Press, 1995.

\bibitem{zbMATH05015420}
S.~{Majid}.
\newblock {What is .. a quantum group?}
\newblock {\em {Notices Am. Math. Soc.}}, 53(1):30--31, 2006.

\bibitem{Moore1}
J.~C. Moore.
\newblock Homotopie des complexes mono\"\i daux, i.
\newblock {\em S\'eminaire Henri Cartan}, 7(2):1--8, 1954-1955.
\newblock talk:18.

\bibitem{zbMATH06520714}
J.~C. {Morton} and R.~{Picken}.
\newblock {Transformation double categories associated to 2-group actions.}
\newblock {\em {Theory Appl. Categ.}}, 30:1429--1468, 2015.

\bibitem{MP1}
A.~{Mutlu} and T.~{Porter}.
\newblock {Iterated Peiffer pairings in the Moore complex of a simplicial
  group.}
\newblock {\em {Appl. Categ. Struct.}}, 9(2):111--130, 2001.

\bibitem{zbMATH03957392}
T.~{Porter}.
\newblock {Extensions, crossed modules and internal categories in categories of
  groups with operations.}
\newblock {\em {Proc. Edinb. Math. Soc., II. Ser.}}, 30:373--381, 1987.

\bibitem{Menagerie}
T.~{Porter}.
\newblock {\em {The Crossed Menagerie: an introduction to crossed gadgetry and
  cohomology in algebra and topology}}.
\newblock \text{Available at:
  \url{http://ncatlab.org/nlab/files/menagerie12a.pdf}}, 2018.

\bibitem{Rad1}
D.~E. {Radford}.
\newblock {The structure of Hopf algebras with a projection.}
\newblock {\em {J. Algebra}}, 92:322--347, 1985.

\bibitem{zbMATH00427768}
D.~E. {Radford} and J.~{Towber}.
\newblock {Yetter-Drinfel'd categories associated to an arbitrary bialgebra.}
\newblock {\em {J. Pure Appl. Algebra}}, 87(3):259--279, 1993.

\bibitem{SW1}
M.~{Sweedler}.
\newblock {\em Hopf algebras}.
\newblock Mathematics lecture note series. W. A. Benjamin, 1969.

\bibitem{zbMATH01584734}
M.~{Takeuchi}.
\newblock {Survey of braided Hopf algebras.}
\newblock In {\em {New trends in Hopf algebra theory. Proceedings of the
  colloquium on quantum groups and Hopf algebras, La Falda, Sierras de
  C\'ordoba, Argentina, August 9--13, 1999}}, pages 301--323. Providence, RI:
  American Mathematical Society (AMS), 2000.

\bibitem{zbMATH06738086}
V.~{Turaev} and A.~{Virelizier}.
\newblock {\em {Monoidal categories and topological field theory.}}, volume
  322.
\newblock Basel: Birkh\"auser/Springer, 2017.

\bibitem{WCH2}
J.~Whitehead.
\newblock Combinatorial homotopy. ii.
\newblock {\em Bull. Amer. Math. Soc.}, 55(5):453--496, 05 1949.

\bibitem{zbMATH04171177}
D.~N. {Yetter}.
\newblock {Quantum groups and representations of monoidal categories.}
\newblock {\em {Math. Proc. Camb. Philos. Soc.}}, 108(2):261--290, 1990.

\end{thebibliography}
%\nocite{*}

\end{document}